\documentclass[a4paper, 10pt, 
]{scrartcl}

\usepackage{amsmath, amsthm, amssymb}  
\usepackage[latin1]{inputenc}
\usepackage[T1]{fontenc}
\usepackage{arydshln}
\usepackage{multirow}
\usepackage{authblk}
\usepackage{graphicx}

\setkomafont{sectioning}{\normalcolor\bfseries}

\renewcommand{\S}{\mathcal{S}^{(i)}}

\newcommand{\halb}{\frac{1}{2}} 
\begin{document}

\begin{center}
	{ \LARGE \bf A square entropy stable flux limiter\\ for $P_NP_M$ schemes}\\

	\hspace{15mm}

	{\large
	Claus R. Goetz
		\footnote{ 
			Corresponding author, 
			{\sc 
				University of Hamburg,  
				Department of Mathematics, 
				Bundesstr. 55, 
				D-20146 Hamburg, Germany,} 
				\texttt{claus.goetz@uni-hamburg.de} 
			} \footnotemark[2]
  	and
  	Michael Dumbser
  		\footnote{
			{\sc
			 	University of Trento, 
		 		Department of Civil, Environmental and Mechanical Engineering, 
				Via Mesiano 77, 
				I-38123 Trento, Italy, 
			}
			\texttt{michael.dumbser@unitn.it}, 
		}
	}	
\end{center}

\hspace{10mm}

\begin{abstract}
We study some theoretical aspects of $P_NP_M$ schemes, which are a novel class of high order accurate reconstruction based discontinuous Galerkin (DG) schemes for hyperbolic conservation laws. The $P_NP_M$ schemes store and evolve the discrete solution $u_h$ under the form of piecewise polynomials of degree $N$, while piecewise polynomials $w_h$ of degree $M \geq N$ are used for the computation of the volume and boundary fluxes. 
The piecewise polynomials $w_h$ are obtained from $u_h$ via a suitable reconstruction or recovery operator. 
The $P_NP_M$ approach contains high order finite volume methods ($N=0$) as well as classical DG schemes ($N=M$) as special cases of a more general framework. Furthermore, for $N \neq M$ and $N>0$, 
it leads to a new intermediate class of methods, which can be denoted either as Hermite finite volume or as reconstructed DG methods. We show analytically why $P_NP_M$ methods for $N \neq M$ are, 
in general, not $L^2$-diminishing. 
To this end, we extend the well-known cell entropy inequality and the following $L^2$ stability result of Jiang and Shu for DG methods (i.e. for $N=M$) to the general $P_NP_M$ case and identify 
which part in the reconstruction step may cause the instability. With this insight we design a flux limiter that enforces a cell square entropy inequality and thus an $L^2$ stability condition 
for $P_NP_M$ schemes for scalar conservation laws in one space dimension. Furthermore, in this paper we prove existence and uniqueness of the solution of the $P_NP_M$ reconstruction operator. 
\end{abstract}

\section{Introduction}

Two of the most popular families of high order schemes for hyperbolic conservation laws are either reconstruction-based \emph{finite volume} schemes, 
or \emph{discontinuous Galerkin (DG)} methods that  directly evolve the expansion coefficients of a piecewise polynomial data representation, and thus 
in principle do not need any reconstruction stage. 
The $P_NP_M$ philosophy introduced in \cite{Dumbser2008,ADERNSE} provides a unified framework for the treatment of both approaches. At this point it is 
important to stress that reconstruction operators have already been used in connection with the DG method before, but only in the context of accuracy-enhancing  
\textit{postprocessors}, which are applied to the discrete solution only once at the end of the simulation, see \cite{ryan,cb-suli-shu}, or as 
\textit{nonlinear moment limiters}, in order to control spurious oscillations at shock waves, see \cite{QiuShu2,QiuShu3,balsara2007}. Instead, the use of 
a reconstruction operator as a mechanism to enhance the accuracy of the DG scheme within each time step has been proposed for the first time within the $P_NP_M$ 
framework \cite{Dumbser2008,ADERNSE}. 
Similar reconstruction-based DG schemes have subsequently been developed and applied also for example in \cite{Luo1,Luo2,HybridDG}. 

In the context of finite volume schemes, high order of accuracy in space is achieved by reconstructing a high order piecewise polynomial representation of the solution (say, an ENO or WENO polynomial) 
from the known cell averages \cite{harteneno,shu_efficient_weno,balserashu,TitarevToroWENO3D}. An example for a high order fully discrete one-step method following this approach is the ADER 
scheme that is based on the approximate solution of generalized Riemann problems arising from the piecewise polynomial reconstruction of the data, see e.g. 
\cite{toro3,titarevtoro,Toro:2006a,AboiyarIske,CastroToro,GoetzIske} for more details. Methods that only evolve cell averages in time have the advantage that they can use large time-steps (CFL $<$ 1), but this comes at the cost of requiring large reconstruction stencils.  

In DG schemes, on the other hand, the high order spatial accuracy derives from the fact that the discrete solution of the PDE is directly represented and evolved terms of piecewise polynomials of degree $N$ in space. At very high order, however, this requires very small time steps (CFL $< 1/(2N+1)$). A natural approach for reaching higher accuracy in the DG framework is then to also include a suitable reconstruction procedure in the scheme. Since more degrees of freedom are available in each cell, we can work with compact reconstruction stencils and still keep a larger CFL number, see \cite{Dumbser2008}. This is the basic idea of the $P_NP_M$ method: The solution is represented in a finite element space of piecewise polynomials of degree $N\geq 0$ (hence the $P_N$ in the name of the method) and at each time step before the time evolution is carried out, a high order reconstruction of piecewise polynomials of degree $M \geq N$ is computed (hence the $P_M$). In this framework, the pure DG method can be viewed as a $P_NP_N$ scheme, while the case $P_0P_M$ corresponds to the high order finite volume schemes. For $N>0$ and $M>N$ a family of hybrid schemes emerges.

The $P_NP_M$ method has already been successfully applied to a variety of complex flow problems, including viscous terms, non-conservative products and stiff source terms, see 
\cite{Dumbser2008,ADERNSE,USFORCE2,DumbserZanotti}, but our 
understanding of its analytical properties has not reached a mature level yet. In particular, it was observed that even in the $1D$ case for scalar, linear problems, 
the method is in general not strictly $L^2$-diminishing \cite{GassnerDipl,dissdumbser} and the famous cell square entropy inequality of Jiang and Shu \cite{JiangShu94} 
for DG schemes is not valid for $M>N$. 

In this paper, we extend the technique of Jiang and Shu to the case $M>N$ and derive a modified sufficient criterion for entropy stability that reduces to the original criterion of Jiang and Shu for $M = N$. We demonstrate why this criterion is in general not satisfied for problems with large jumps in the solution. In the linear case the new condition is a simple algebraic relation between the jump in the reconstruction and the jump in  the data. Those jumps can be expected to be of similar size for smooth solutions, but can differ substantially in the vicinity of strong shocks. In order to stabilize the method, we propose a new \emph{flux limiter} that strictly enforces the cell entropy inequality and thus $L^2$ stability of general $P_NP_M$ schemes. 

The rest of this paper is organized as follows: At first we describe the $P_NP_M$ method in  Section \ref{sec:PNPM}.  After that we discuss the reconstruction procedure in Section \ref{sec:Rec}. In section \ref{sec:Entropy}, we derive a condition on the reconstruction operator in order to satisfy a cell entropy inequality. Our approach is analogous to the technique of Jiang and Shu and we discuss which steps in their proof have to be modified when we include a reconstruction operator in the scheme. We use the new insight where entropy stability may be violated to develop a new flux limiter that ensures entropy stability in Section \ref{sec:limiter}. The new limiter is tested and validated in several numerical examples, which are presented in Section \ref{sec:numerics}. We summarize the work and draw conclusion in Section \ref{sec:conclusion}.

Finally, in the appendix, we provide a proof that the $P_NP_M$ reconstruction problem on the three cell central stencil has a unique solution for $M=3N+2$, which is the maximal possible $M$ that can be reached in this case. It seems that such a result for arbitrary $N$ has not been available in the literature so far.

\section{The $P_NP_M$ method} \label{sec:PNPM}
In this section we describe the fundamental concepts of the $P_NP_M$ method. Our focus lies on the analytical aspects of the scheme, so questions of practical implementation will not be addressed in detail. For those we refer to \cite{Dumbser2008}. 
We consider the following framework: 

We solve the initial value problem for a scalar conservation law in one space dimension,
\begin{align}
\label{CL}
&	\frac{\partial}{\partial t}u(x,t) + \frac{\partial}{\partial x}f \left( u(x,t) \right) = 0, \qquad (x,t) \in \Omega \times [0,\infty),\\
&	\quad \Omega \subseteq \mathbb{R}, \qquad u(x,t) \in \mathcal{U} \subset \mathbb{R}, \qquad f: \mathcal{U} \to \mathbb{R} , \notag
\end{align}
with the space of admissible states $\mathcal{U}$, initial data 
\[
	u(x,0) = \bar{u}(x)
\]
and, if necessary, suitable boundary conditions.
Let $\mathcal{T} = \{ T^{(i)}, ~i \in I\}$ for some index set $I$ be a partition of $\Omega$. In the following we assume for simplicity that we have a uniform grid on the whole real line,
\[
	\Omega = \mathbb{R}, \quad T^{(i)} = [x_{i-\halb}, x_{i+\halb}], \quad \text{ with } x_{i+\halb} = \left( i + \halb \right) h, ~ i\in \mathbb{Z},~ h>0.
\]
Let $u_h(x,t) \in \mathcal{V}_h$ be a piecewise polynomial representation of the solution at time $t$, such that in each cell $T^{(i)} = [ x_{i-\halb}, x_{i+\halb}]$ 
the function $u^{(i)} = \left. u_h \right|_{T^{(i)}}$ is a polynomial of degree $N$. At time $t=0$, $u_h(x,0)$ is obtained as the $L^2$ projection of the initial data $u(x,0)$ 
onto the space $\mathcal{V}_h$ of piecewise polynomials of degree $N$ on $\mathcal{T}$.   

We denote one-sided limits at the cell interfaces by
\[
	u_{i-\halb}^+(t) = \lim_{\substack{ x \to x_{i-\halb} \\ x > x_{i-\halb} }} u_h(x,t), \qquad
	u_{i+\halb}^-(t) = \lim_{\substack{ x \to x_{i+\halb} \\ x < x_{i+\halb} }} u_h(x,t).
\]
Moreover, let $w_h(x,t) \in \mathcal{W}_h$ be a function that is reconstructed from $u_h(x,t)$ on a reconstruction stencil $\mathcal{S}^{(i)} = \{T^{(i-k)}, \dots, T^{(i+\ell)}\}$, 
with the help of a reconstruction operator $\mathcal{R}$ (to be specified later), 
\[
	w^{(i)} = \mathcal{R}(u^{(i-k)},\dots, u^{(i+\ell)}), \qquad k,\ell \geq 0,
\] 
with $w^{(i)} = \left. w_h \right|_{T^{(i)}}$, so that $w_h$ is given in each $T^{(i)}$ by a polynomial of degree $M \geq N$.  With $\mathcal{W}_h$ we denote the space of piecewise polynomials 
of degree $M$ on $\mathcal{T}$. In this paper we focus on the \emph{three point cell central stencil}, 
i.e. $\mathcal{S}^{(i)}=\{ T^{(i-1)}, T^{(i)}, T^{(i+1)}\}$. The reconstruction method according to  \cite{Dumbser2008} is described in the following Section \ref{sec:Rec}, but a fundamental 
requirement of the reconstruction operator is a generalized conservation principle in the sense that  
\begin{equation}
\label{eqn.recprop} 
   \int \limits_{T^{(i)}} \varphi(x) w_h(x,t) \, dx = \int \limits_{T^{(i)}} \varphi(x) u_h(x,t) \, dx, \qquad \forall \varphi(x) \in \mathcal{V}_h.   
\end{equation} 

The (semi-discrete) $P_NP_M$ scheme is then derived from a weak formulation of the conservation law \eqref{CL} by integrating over a control volume 
$T^{(i)}$ after multiplication with a test function $\varphi=\varphi(x)$ from a suitable space and after integration by parts of the term containing the flux derivative: 
\begin{align}\label{weak 1}
	\int \limits_{T^{(i)}} \frac{\partial w_h}{\partial t} \varphi(x)~dx + f_{i+\halb}^w(t)\varphi_{i+\halb}^-  - f_{i-\halb}^w(t) \varphi_{i-\halb}^+  
	- \int \limits_{T^{(i)}} \!\! f \left( w_h(x,t) \right) \frac{d \varphi}{dx} ~dx = 0.
\end{align}
Here, $f_{i+\halb}^w(t)$ is given by 
\[
	f_{i+\halb}^w(t) =  \bar{f}_{i+1/2}\left(w_{i+\halb}^-(t),~ w_{i+\halb}^+(t) \right),	
\]
where $\bar{f}_{i+\halb}$ is a Lipschitz continuous and consistent numerical flux, i.e. $\bar{f}_{i+\halb}(u,u) = f(u)$. The superscript $w$ indicates that the numerical flux is evaluated at reconstructed values, i.e. using the boundary extrapolated values of the reconstructed solution $w_h$ at the cell interfaces. 
Moreover, suppose that $\bar{f}_{i+\halb}$ is an \emph{E-flux} in the sense of Osher \cite{oshereflux}:
\[
	\left( \bar{f}_{i+\halb}(u^-,~u^+) - f(\xi) \right)(u^+ - u^-) \leq 0 \qquad \text{ for all } \xi \text{ between } u^- \text{ and } u^+.
\]
Finally, choose a basis $\{ \Phi_\ell^{(i)}~:~\ell = 0, \dots, N\}$ for the space of polynomials of degree $N$ on $T^{(i)}$ and write
\[
	u^{(i)} = \sum_{\ell = 0}^N \hat{u}_\ell^{(i)}(t) \Phi_\ell^{(i)}(x).
\]
For $\varphi \in \mathcal{V}_h$, the property of the reconstruction operator \eqref{eqn.recprop} can be applied and thus 
the weak formulation \eqref{weak 1} becomes 
\begin{align}\label{weak_1b}
	\int \limits_{T^{(i)}} \frac{\partial u_h}{\partial t} \varphi(x) ~ dx + f_{i+\halb}^w(t)\varphi_{i+\halb}^- - f_{i-\halb}^w(t) \varphi_{i-\halb}^+  
	- \int \limits_{T^{(i)}} \!\! f \left( w_h(x,t) \right) \frac{d \varphi}{dx}  ~dx = 0, \nonumber \\ \qquad \forall \varphi \in \mathcal{V}_h. 
\end{align}
Let $\varphi = \Phi_k^{(i)}$ for $k=0,\dots,N$ and thus \eqref{weak_1b} leads to an ODE system for the temporal evolution of the degrees of freedom $\hat{u}_\ell^{(i)}(t)$. 
Note that \eqref{weak 1} and \eqref{weak_1b} belong to the class of Petrov-Galerkin schemes, were the discrete solution $w_h \in \mathcal{W}_h$ and 
the test functions $\varphi \in \mathcal{V}_h$ are from different spaces. 

\section{Reconstruction} \label{sec:Rec}

Our goal in this section is to construct a function $w_h$, such that for each cell $T^{(i)}$ the function $w^{(i)} = \left. w_h \right|_{T^{(i)}}$ is a polynomial of degree $M>N$. To this end, we first choose a \emph{reconstruction stencil} $\mathcal{S}^{(i)} \subset \mathcal{T}$ for each $T^{(i)}$ and construct a polynomial $w_s^{(i)}$ of degree $M$ \emph{on the whole stencil} $\mathcal{S}^{(i)}$. 
The function $w^{(i)} = w_s^{(i)}|_{T^{(i)}}$ is then given as the restriction of that reconstruction polynomial $w_s^{(i)}$ to the cell $T^{(i)}$, while $w_h$ is given by the union of all the $w^{(i)}$ for all $T^{(i)} \in \mathcal{T}$. For simplicity, we only discuss the reconstruction procedure on the three cell central stencil $\mathcal{S}^{(i)} = \{ T^{(i-1)},~ T^{(i)},~ T^{(i+1)}\}$. Here, by a slight abuse of notation, we also write $\mathcal{S}^{(i)} = T^{(i-1)} \cup T^{(i)} \cup T^{(i+1)}$ to denote the union of all cells in the stencil. Denote the space of piecewise polynomials of degree $N$ on $\mathcal{S}^{(i)}$ by
\[
	V_h^N\left( \mathcal{S}^{(i)} \right) = \left\{ p: \mathcal{S}^{(i)} \to \mathbb{R} ~ : ~ \left. p\right|_{ T^{(i)} } \in \mathbb{P}^N \right\} ,
\]
and the space of polynomials of degree $M>N$ on $\mathcal{S}^{(i)}$ by
\[
	W_h^M\left( \mathcal{S}^{(i)} \right) = \left\{ p: \mathcal{S}^{(i)} \to \mathbb{R} ~ : ~  p \in \mathbb{P}^{M} \right\}.
\]
Note that while in the end we are looking for a reconstruction in the space of \emph{piecewise} polynomials of degree $M$ (i.e. a function in $\mathcal{W}_h = V_h^M(\mathcal{T})$), for the moment we consider the space $W^{M}\left( \mathcal{S}^{(i)} \right)$ of polynomials that are continuously differentiable on the whole stencil $\mathcal{S}^{(i)}$.

For each cell $T^{(j)}\in\mathcal{S}^{(i)}$, let $\{ \Phi_\ell^{(j)},~\ell \in \{0, \dots, N \} \}$ be the shifted Legendre polynomials on $T^{(i)}$, which form a basis of the polynomial space 
$\mathbb{P}^N\left( T^{(j)} \right)$, that is orthogonal in the local $L^2$-sense:
\[
	\left\langle \Phi_m^{(j)}, \Phi_n^{(j)} \right\rangle_{T^{(j)}} = 
	\int_{T^{(j)} }	\Phi_m^{(j)}(x) \Phi_n^{(j)}(x)~dx = 0, 
	\qquad \text{ if } n \neq m.
\]

For $W_h^{M}\left( \mathcal{S}^{(i)} \right)$ we choose a basis $\left\{ \Psi_k ,~ k \in \{0, \dots, M\} \right\}$ such that
\[
	\left\langle \Psi_m, \Psi_n \right\rangle_{T^{(i)}} = 0 \qquad \text{ for } n \neq m.
\]
Note that while all $\Psi_k$ are defined on the whole stencil $\mathcal{S}^{(i)}$, we require only their restrictions to the central cell $T^{(i)}$ to be orthogonal. To achieve this, we simply take the $\Psi_k$ to be Legendre polynomials on the central cell $T^{(i)}$ and extend them to the whole stencil. Note that this in particular means that
\[
	\left.\Psi_k \right|_{T^{(i)}} = \Phi_k^{(i)}, \qquad \text{ for } k = 0, \dots, N. 
\]

In the following, we denote by $w_s^{(i)}$ the continuous extension of the polynomial $w^{(i)}$ of degree $M$ to the entire stencil $\mathcal{S}^{(i)}$, while 
$\left. w_s^{(i)} \right|_{T^{(i)}} = w^{(i)}$.  
According to \cite{Dumbser2008}, we require the reconstruction $w_s^{(i)} \in W^{M}\left( \mathcal{S}^{(i)} \right)$ from a given $u_h \in V_h^N\left( \mathcal{S}^{(i)} \right)$, to satisfy a 
\emph{generalized conservation principle} or \textit{identity in the weak sense}, i.e. we require \eqref{eqn.recprop} to hold for all cells in $\mathcal{S}^{(i)}$. That means 
\[
	\left\langle w_s^{(i)}, \Phi_\ell^{(j)} \right\rangle_{T^{(j)}} = \left\langle u_h, \Phi_\ell^{(j)} \right\rangle_{T^{(j)}}
	\qquad \text{ for each } \ell \in \{0, \dots, N\}, ~ j \in \{i-1, i, i+1\},
\]
With 
\[
	w_s^{(i)} = \sum_{k=0}^{M} \hat{w}_k^{(i)} \Psi_k, \qquad u^{(j)} = \sum_{k=0}^{N} \hat{u}_k^{(j)} \Phi_k^{(j)},
\]
for each $\ell \in \{0, \dots, N\}$ and each $j\in \{i-1,i,i+1\}$ this leads to the condition
\begin{equation}
	\sum_{k=0}^{M} \hat{w}^{(i)}_k \langle \Psi_k, \Phi_\ell^{(j)} \rangle_{T^{(j)}} = \sum_{k=0}^{N} \hat{u}_k^{(j)} \langle  \Phi_k^{(j)},  \Phi_\ell^{(j)} \rangle_{T^{(j)}}.
	\label{inverse L2 - 1}
\end{equation}

Let's define a basis of $V_h^N\left( \mathcal{S}^{(i)}\right)$ in the following way: At first, we extend each of the local basis functions $\Phi_\ell^{(j)}$ to the whole stencil $\mathcal{S}^{(i)}$ by setting
\[
	\tilde{\Phi}_\ell^{(j)}(x) =~ \left\{
	\begin{array}{cc}
		\Phi_\ell^{(j)}(x), & \quad x \in T^{(j)},\\
		0, 					& \quad \text{ else }.
	\end{array}
	\right.
\]
Next, let $k$ be an index $k \equiv k(\ell, j)$, such that
\begin{align*}
	k(0,i-1) & = 0,    \quad	\dots, 	& k(N, i-1)	& = N,  &\\ 
	k(0, i ) & = N+1,  \quad 	\dots, 	& k(N, i )	& = 2N+1, &\\
	k(0,i+1) & = 2N+2, \quad	\dots, 	& k(N, i+1) & = 3N+2. &
\end{align*}
Finally, we set 
\[
	\Theta_k = \tilde{\Phi}_\ell^{(j)}, \qquad \ell \in \{0, \dots, N\}, ~ j \in \{i-1, i, i+1\},
	\qquad  k = k(\ell,j).
\]
Then the $\Theta_k$ form a basis of $V_h^N\left( \mathcal{S}^{(i)} \right)$, that is orthogonal in the sense that
\[
	\left\langle \Theta_m, \Theta_n \right\rangle_{\mathcal{S}^{(i)}} = 
	\int_{\mathcal{S}^{(i)} } \Theta_m(x) \Theta_n(x)~dx = 0, 
	\qquad \text{ if } m \neq n.
\]
With this, we can rewrite \eqref{inverse L2 - 1} as
\begin{equation}
	\sum_{\ell=0}^{M} \hat{w}_\ell^{(i)} \langle \Psi_\ell, \Theta_k \rangle_{\mathcal{S}^{(i)}} = \sum_{m=0}^{3N+2} \hat{u}_m \langle  \Theta_m,  \Theta_k \rangle_{\mathcal{S}^{(i)}},
	\qquad k \in \{0,\dots,3N+2\}.
	\label{inverse L2 - 2}
\end{equation}

If $M < 3N+2$, this is an overdetermined system, which we solve at the aid of a constrained least-squares method, see \cite{DumbserKaeser06b, Dumbser2008}. In Appendix A we give an existence and uniqueness proof for $M=3N+2$, which, to our knowledge, has not been available in the literature so far. 

In particular, by requiring conservation on the central cell of the stencil, we get
\[
	\hat{w}_\ell^{(i)} = \hat{u}_\ell^{(i)}, \qquad \text{ for } \ell=0,\dots,N,
\]
which in turn implies that the \emph{residual} $r_h = w_h - u_h$, is orthogonal to $\mathcal{V}_h$ and thus to $u_h$ on $T^{(i)}$:
\[
\label{eqn.ortho} 
	\int_{T^{(i)}} r_h(x) u_h(x)~dx = \int_{T^{(i)}} \left( \sum_{k=N+1}^M \hat{w}_k^{(i)} \Psi_k(x) \right) \left( \sum_{k=0}^N \hat{u}_k^{(i)} \Phi_k^{(i)}(x) \right)~dx = 0.
\]

\section{The cell entropy condition}\label{sec:Entropy}

The focus of our analysis in this paper is the \emph{entropy stability} (or lack thereof) of $P_NP_M$ schemes. Recall that a convex function $Q: \mathcal{U}\to\mathbb{R}$ is called an \emph{entropy} for the conservation law \eqref{CL}, if there exists a so-called \emph{entropy flux} $F$ with $F'(u) = Q'(u)f'(u)$. We say that a weak solution $u$ of \eqref{CL} satisfies the \emph{entropy inequality}, if
\begin{equation}
\label{entropy inequality}
	\frac{\partial }{\partial t}Q(u) + \frac{\partial}{\partial x}F(u) \leq 0
\end{equation}
in a distributional sense. For classical (i.e. $C^1$) solutions, \eqref{entropy inequality} is satisfied with equality. A discrete version of \eqref{entropy inequality} is given for each cell $T^{(i)}\in \mathcal{T}$ by
\begin{equation}
\label{num entropy inequality}
\int_{T^{(i)}} \frac{\partial}{\partial t}Q\left( u_h(x,t) \right)~dx + F_{i+\halb}(t) - F_{i-\halb}(t) \leq 0,
\end{equation}
where $F_{i+\halb}$ is a numerical entropy flux, consistent with $F$. For the special case of the \emph{square entropy} $Q(u) = u^2/2$, by summing over all $T^{(i)} \in \mathcal{T}$, condition \eqref{num entropy inequality} leads to the $L^2$-stability of the scheme:
\[
	\frac{d}{dt} \| u_h(x,t) \|_{L^2(\Omega)} \leq 0.
\]
Although we will mostly talk about entropy stability in this paper, it is in fact the $L^2$-stability of the scheme that we are most interested in.

Jiang and Shu \cite{JiangShu94} have shown that the pure discontinuous Galerkin scheme, which we interpret as a $P_NP_N$ scheme with $M=N$, i.e. with $w_h=u_h$ in \eqref{weak_1b}, satisfies the 
discrete cell entropy inequality \eqref{entropy inequality} for the square entropy. This does, however, not hold in general for $w_h \neq u_h$. To our knowledge, this phenomenon has so far not yet been explained  theoretically. 
In this section we follow the original construction of Jiang and Shu and highlight which steps need to be modified to cover the general $P_NP_M$ case. We show why the reconstruction procedure 
can lead to a violation of the entropy condition. 

We begin analogously to Jiang and Shu, i.e. in \eqref{weak_1b} take $\varphi = u_h(x,t)$, so that one gets  
\begin{multline}\label{weak 2}
	\frac{d}{dt}\int \limits_{T^{(i)}} \frac{ \left( u_h(x,t) \right)^2}{2} ~dx + f_{i+\halb}^w(t)u_{i+\halb}^-(t) - f_{i-\halb}^w(t)u_{i-\halb}^+(t)  
	- \int \limits_{T^{(i)}} f \left( w_h(x,t) \right) \frac{\partial u_h }{\partial x }~dx = 0.
\end{multline}
Let us denote
\[
	E^{(i)}(t) = \int_{T^{(i)}}\frac{\left( u_h(x,t)\right)^2}{2}~dx.
\]
Recall that in the $P_NP_M$ method, contrary to the pure DG method, numerical fluxes and the physical flux in the volume integral are evaluated at the reconstructed polynomial $w_h$, 
rather than at the DG polynomial $u_h$. In general $u_h \in \mathcal{V}_h$ and $w_h \in \mathcal{W}_h$ are not from the same space. However, we formally interpret $f_{i+\halb}^w$ 
as a four-point flux $f_{i+\halb}^{\mathcal{R}(u)}$: 
\begin{align*}
	f_{i+\halb}^{\mathcal{R}(u)}(t) 
&	= f_{i+\halb}^{\mathcal{R}(u)}\left(u^{(i-1)}, \dots, u^{(i+2)} \right)(t)\\ 
&	= \bar{f}_{i+\halb}\left( \mathcal{R}\left(u^{(i-1)}, u^{(i)}, u^{(i+1)}\right)_{i+\halb}^-(t), ~\mathcal{R}\left(u^{(i)}, u^{(i+1)}, u^{(i+2)}\right)_{i+\halb}^+(t) \right). 
\end{align*}
So we can formally do all analysis in the space of polynomials of degree $N$. The consistency condition
\[
	f_{i+\halb}^{\mathcal{R}(u)}(u, \dots,u) = f(u)
\]
is satisfied if the reconstruction operator preserves constants, i.e.
\[
	\mathcal{R}(u,\dots,u) = u.
\]

The next step is to construct a numerical entropy flux. Recall that an entropy flux $F$ for the square entropy is given by
\[
	F(u) = f(u)u-g(u)
\]
where $g$ is a primitive of $f$:
\[
	g(u) = \int f(u)~du.
\]
To simplify the notation, in the following let us drop the time-dependence in the notation, but keep in mind that we are describing a \emph{semi-discrete} method.

We are looking for a numerical entropy flux $F_{i+\halb}$ that is consistent with $F$. Take $g$ as defined above and write $f(w) = f(u) + (f(w)-f(u))$. Then we get from \eqref{weak 2}:
\begin{multline}\label{weak 3}
	\frac{d}{dt}E^{(i)}	+ f_{i+\halb}^wu_{i+\halb}^- 
	- f_{i-\halb}^wu_{i-\halb}^+  + g(u_{i-\halb}^+) - g(u_{i+\halb}^-) 
	- \int \limits_{T^{(i)}} \! \left( f\left( w_h \right) - f\left( u_h \right) \right) \frac{\partial u_h}{\partial x } ~dx = 0.
\end{multline}
Adding and subtracting $f_{i-\halb}^w u_{i-\halb}^- - g(u_{i-\halb}^-)$ from \eqref{weak 3} and rearranging terms yields 
\begin{multline}\label{weak 3b}
	\frac{d}{dt}E^{(i)}	+ f_{i+\halb}^wu_{i+\halb}^- - g(u_{i+\halb}^-) - f_{i-\halb}^wu_{i-\halb}^-   + g(u_{i-\halb}^-)  
	- \int \limits_{T^{(i)}} \! \left( f\left( w_h \right) - f\left( u_h \right) \right) \frac{\partial u_h}{\partial x } ~dx, \\
		- f_{i-\halb}^wu_{i-\halb}^+  + f_{i-\halb}^wu_{i-\halb}^-   - g(u_{i-\halb}^-)   + g(u_{i-\halb}^+)  = 0. \qquad \qquad  
\end{multline}
Setting
\[
	\tilde{F}_{i+\halb} = f_{i+\halb}^wu_{i+\halb}^- - g(u_{i+\halb}^-)
\]
and
\begin{align*}
	V^{(i)} 
&	= \int \limits_{T^{(i)}} \! \left( f\left( w_h \right) - f\left( u_h \right) \right) \frac{\partial u_h}{\partial x } ~dx 
		\\ 
	B^{(i)} 
&	=	- f_{i-\halb}^w \left(u_{i-\halb}^+ - u_{i-\halb}^-\right) + g(u_{i-\halb}^+) - g(u_{i-\halb}^-)- V^{(i)},
\end{align*}
eqn. \eqref{weak 3b} becomes
\begin{align}\label{weak 4}
	\frac{d}{dt}E^{(i)}	+ \tilde{F}_{i+\halb} - \tilde{F}_{i-\halb} + B^{(i)} = 0.
\end{align}
The numerical entropy flux $\tilde{F}_{i+\halb}$ is consistent with the entropy flux $F$. Therefore, if $B^{(i)} \geq 0$, we get the discrete entropy inequality
\begin{align}\label{entropy condition}
	\frac{d}{dt}E^{(i)} + \tilde{F}_{i+\halb} - \tilde{F}_{i-\halb} \leq 0.
\end{align}
The rest of this section is devoted to analysing under which conditions we can guarantee $B^{(i)} \geq 0$ and thereby the entropy inequality \eqref{entropy inequality}.

By the mean-value theorem, there exists a $\xi$ between $u_{i-\halb}^-$ and $u_{i-\halb}^+$, such that
\[
	g(u_{i-\halb}^+) - g(u_{i-\halb}^-) = \int_{u_{i-\halb}^-}^{u_{i-\halb}^+} f(u)~du 
	= \left( u_{i-\halb}^+ - u_{i-\halb}^- \right) f(\xi),
\]
and so a sufficient condition for $B^{(i)} \geq 0$ is
\begin{equation} \label{positive B}
	- \left( u_{i-\halb}^+ - u_{i-\halb}^- \right) \left( f^w_{i-\halb} - f(\xi) \right) 
	- V^{(i)}\geq 0.
\end{equation}
Note that up to here there our analysis is completely analogous to the Jiang and Shu proof \cite{JiangShu94}. The only differences are that in our case 
$f_{i-\halb}^w$ is evaluated using the reconstructed solution $w_h$ and that we have to include the term $V^{(i)}$.

Assume that $u_{i-\halb}^- > u_{i-\halb}^+$, then condition \eqref{positive B} becomes
\begin{align}
	\bar{f}_{i-\halb}\left(w_{i-\halb}^-,~ w_{i-\halb}^+ \right)  \geq 
	\frac{1}{u_{i-\halb}^- - u_{i-\halb}^+} \left( \int_{u_{i-\halb}^+}^{ u_{i-\halb}^-} f(u)~du 
	-V^{(i)}\right).
	\label{entropy condition 2}
\end{align}  
Obviously, this condition depends on the flux $f$, the numerical flux $\bar{f}_{i+\halb}$ and the reconstruction $w_h = \mathcal{R}(u_h)$. In the following, we discuss the restrictions on $w_h$ imposed by \eqref{entropy condition 2} for \emph{given $f$ and $\bar{f}_{i+\halb}$}.

To better understand where the entropy condition may not be satisfied, we start with the easiest case possible: We discuss the linear advection equation, $f(u) = \lambda u,~\lambda>0$ constant and take $\bar{f}_{i+\halb}$ to be the upwind flux. This allows us to find a simple algebraic condition that relates the jump in the reconstruction to the jump in the data. It turns out that already in this simple case, the reconstruction procedure can lead to  a violation of the entropy condition.

Assume that at the cell-interface $i-\halb$ we have a jump in the data with $u_{i-\halb}^- > u_{i-\halb}^+$. With the upwind flux, condition \eqref{entropy condition 2} reads
\[
	\lambda w_{i-\halb}^-  \geq 
	\frac{1}{u_{i-\halb}^- - u_{i-\halb}^+} \left( \int_{u_{i-\halb}^+}^{ u_{i-\halb}^-} f(u)~du 
	-V^{(i)}\right).	
\] 
In the linear case, the flux-average $f(\xi)$ in the jump is simply
\[
	f(\xi) = \frac{1}{u_{i-\halb}^+ - u_{i-\halb}^-} \int_{u_{i-\halb}^-}^{u_{i-\halb}^+}\lambda u~du = \frac{1}{2}\lambda\left( u_{i-\halb}^- + u_{i-\halb}^+ \right),
\]
and the volume integral of the flux difference can be written in terms of the residual $r_h = w_h - u_h$ as 
\begin{equation}
\label{linear volume intreal}
	V^{(i)} =
	 \int_{T^{(i)}} \lambda \left( w_h - u_h \right) \frac{\partial u_h}{\partial x }~dx = \int_{T^{(i)}} \lambda \, r_h \frac{\partial u_h}{\partial x }~dx = 0.
\end{equation}
Since inside each cell $T^{(i)}$ the function $\frac{\partial u_h}{\partial x }$ is a polynomial of degree $k-1$ and therefore can be expressed as a linear combination of  
$\Phi_0^{(i)},\dots,\Phi_{k-1}^{(i)}$, i.e. we have $\frac{\partial u_h}{\partial x } \in \mathcal{V}_h$ and thus the integral \eqref{linear volume intreal} vanishes according 
to the orthogonality of $r_h$ with respect to $\mathcal{V}_h$, see also \eqref{eqn.ortho}. 
In summary, for linear problems \eqref{positive B} reduces to
\[
	\lambda w_{i-\halb}^- - f(\xi) \geq 0.
\]
Thus, we get the \emph{pointwise} condition
\begin{equation}
	w_{i-\halb}^- \geq \frac{1}{2}\left( u_{i-\halb}^- + u_{i-\halb}^+ \right).
	\label{pointwise condition}
\end{equation}
Analogously, if $u_{i-\halb}^- < u_{i-\halb}^+$, we get
\begin{equation}
	w_{i-\halb}^- \leq \frac{1}{2}\left( u_{i-\halb}^- + u_{i-\halb}^+ \right).
		\label{pointwise condition 2}
\end{equation}
It is easy to construct counterexamples, in which \eqref{pointwise condition}, \eqref{pointwise condition 2} is not satisfied. However, numerical experiments suggest, that this occurs in practice only 
when there is a large jump in $u_h$. 

\begin{figure}
\begin{center}
\includegraphics[width=0.6\textwidth]{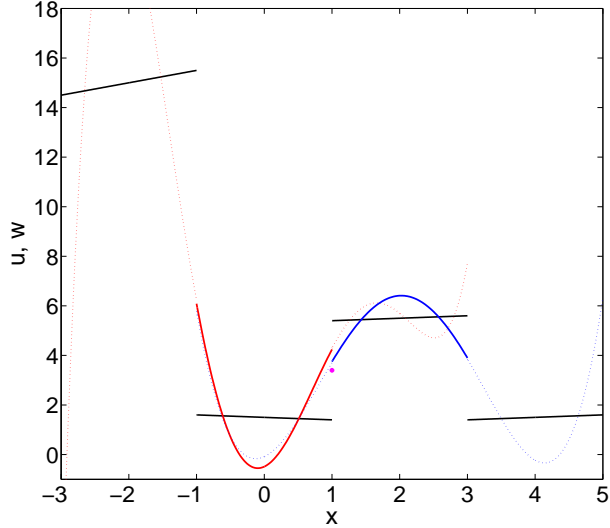}
\end{center}
\caption{A $P_1P_5$ reconstruction that does not satisfy condition \eqref{pointwise condition}, \eqref{pointwise condition 2}. \label{fig:counter example} }
\end{figure}

Figure \ref{fig:counter example} shows an example where condition \eqref{pointwise condition}, \eqref{pointwise condition 2} is not satisfied. We have piecewise linear data (solid black lines) and compute two polynomials of degree five (red and blue). The leftmost there parts of the initial data are used to compute a polynomial $w_L$ (plotted in red), the parts to the right are used for the polynomial $w_R$ (plotted in blue). For the intercell flux between the cell $(-1,1)$ and $(1,3)$ we consider the interface $i - 1/2 = 1$. Condition \eqref{pointwise condition 2} states that in order to satisfy the entropy inequality we need $w_L(1) < \frac{1}{2}(u_1^+ + u_1^-)$, which is marked as a magenta point in the plot. Clearly, this condition is not satisfied. 

In the next section, we describe a way to stabilize the method using a \emph{flux limiter} approach. Our ansatz also covers the nonlinear case.

\section{Flux limiting}\label{sec:limiter}

Now that we have seen where instabilities may occur, we develop a nonlinear fix for that problem. In order to enforce a cell entropy condition, we employ a \emph{flux limiting} approach. Denote by $f_{i+\halb}^u,~f_{i+\halb}^w$ numerical fluxes evaluated using $u_h$ and $w_h$, respectively. Our goal is to find a limiter $\theta_{i+\halb} \in [0,1]$ such that the scheme with the limited flux 
\begin{equation}
	f_{i+\halb}   = f_{i+\halb}^u + \theta_{i+\halb}f_{i+\halb}^r,  \qquad
	f_{i+\halb}^r = f_{i+\halb}^w-f_{i+\halb}^u
	\label{limited flux}
\end{equation}
satisfies the entropy condition. We denote the terms that only include $u_h$ by 
\begin{align*}
	F_{i+\halb}^u 
&	= f_{i+\halb}^u u_{i+\halb}^- - g(u_{i+\halb}^-),\\
	A^{(i)} =
&	- \left(u_{i-\halb}^+ - u_{i-\halb}^-\right)f_{i-\halb}^u + g(u_{i-\halb}^+) - g(u_{i-\halb}^-).
\end{align*}
Note that these are exactly the terms that occur in the original work of Jiang and Shu \cite{JiangShu94} and we recall that the essential point of their proof is that $A^{(i)} \geq 0$ since 
$f_{i-\halb}^u$ is a monotone E-flux. Moreover, let us use the standard notation for a jump at the cell interface 
$	\left[ \left[ u \right]\right]_{i-\halb} = u_{i-\halb}^+-u_{i-\halb}^-$. 
Inserting \eqref{limited flux} in \eqref{weak 4} and rearranging terms yields:
\begin{align}
	\frac{d E^{(i)}}{dt}
&	+ \left( F_{i+\halb}^u + \theta_{i+\halb}f_{i+\halb}^r u_{i+\halb}^- \right)
	- \left( F_{i-\halb}^u + \theta_{i-\halb}f_{i-\halb}^r u_{i-\halb}^- \right) \nonumber \\
& 	+ A^{(i)} - V^{(i)} - \theta_{i-\halb} \left[ \left[ u \right]\right]_{i-\halb} f_{i-\halb}^r  = 0 
\end{align}
Thus, defining a numerical entropy flux as
\[
	\hat{F}_{i+\halb} =  F_{i+\halb}^u + \theta_{i+\halb}f_{i+\halb}^r u_{i+\halb}^-,
\]
the entropy condition
\[
	\frac{d E^{(i)}}{dt} + \hat{F}_{i+\halb} - \hat{F}_{i-\halb} \leq 0
\]
is satisfied if
\begin{equation}
	A^{(i)} - V^{(i)} - \theta_{i-\halb} \left[ \left[ u \right]\right]_{i-\halb} f_{i-\halb}^r  \geq 0. 
\label{nonlinear limiter condition}
\end{equation}
In the linear case $f(u) = \lambda u$ we have $V^{(i)}=0$ due to the orthogonality of the residual $r_h$ w.r.t. $\mathcal{V}_h$, hence 
the condition \eqref{nonlinear limiter condition} simplifies to .
\begin{equation}
	A^{(i)} - \theta_{i-\halb} \left[ \left[ u \right]\right]_{i-\halb} f_{i-\halb}^r  \geq 0. 
\label{linear limiter condition}
\end{equation}
Note that $A^{(i)} \geq 0$ and so in the linear case we can always find a $\theta_{i-\halb} \in [0,1]$ that satisfies \eqref{linear limiter condition}. 

In the nonlinear case, we have to account for the volume integral $V^{(i)}$ of the flux difference in \eqref{nonlinear limiter condition}. 
Since the limiter only acts on the \emph{cell-boundary}, it is possible to include this new term, which only takes information from \emph{inside the cell}. However, it may happen that condition \eqref{nonlinear limiter condition} cannot be satisfied with a $\theta_{i-\halb} \in [0,1]$. In this case, we found that good numerical results can be obtained by setting 
\[
	\sigma_{i-\halb} 
 	= \frac{A^{(i)} - V^{(i)} }{\left[ \left[ u \right]\right]_{i-\halb} f_{i-\halb}^r }
\]
(avoiding divisions by zero) and then compute the actual flux limiter from 
\[ 	
	\theta_{i-\halb} 
	= \max \left( \min ( \sigma_{i-\halb}, 1), 0 \right).
\]
In the update for the degrees of freedom inside each cell, we do not limit the polynomial $w_h$ inside the cell if $\theta_{i-\halb} \in (0,1]$. However, for $\theta_{i-\halb} = 0$, it can happen that $A^{(i)}-V^{(i)}<0$ and in this case the flux limiter alone does not guarantee the entropy condition. If this occurs, additional limiting of the polynomial inside the cell is necessary. Numerical experiments  suggest that in the case of entropy violation for $\theta_{i-\halb}=0$, one can simply set $ \left. w_h \right|_{T^{(i)}} := u^{(i)} + \theta_i (w^{(i)} - u^{(i)}) $ with $\theta_i \in [0,1]$, so that 
$A^{(i)} + V^{(i)} \geq 0$ and thus entropy stability is guaranteed. Numerical evidence shows that even setting $\theta_i = 0$ does not hurt the accuracy too much. 
Just as in the pure DG case, it is important to stress that even with the limiter developed above, we can only enforce the entropy inequality and thereby the $L^2$-stability of the scheme. It does not take into account requirements on the $L^\infty$ norm. Moreover, the reconstruction method we use is linear. Due to Godunov's theorem \cite{godunov}, the high order $P_NP_M$ schemes with this reconstruction cannot be monotone. Additional limiting or a different, nonlinear WENO/HWENO reconstruction procedure \cite{QiuShu2,QiuShu3,balsara2007} are necessary to deal with spurious 
oscillations at shock waves. 

\section{Numerical results}\label{sec:numerics}
\subsection{Linear advection with smooth initial data}
At first we test the effects of the flux limiter for a linear problem with smooth initial data. In this case, even the unlimited $P_NP_M$ scheme has little to no problems with $L^2$ stability and so the main question will be whether the flux limiting affects the order of accuracy. Consider the problem
\[
	\frac{\partial u}{\partial t} + \frac{\partial u}{\partial x} = 0, \qquad x\in (-1,1)
\]
with initial data
\[
	u(x,0) = \sin (\pi x)^4
\]
and periodic boundary conditions. Errors are measured at $t=1$. Time integration is performed with an $(M+1)$-stage linear Runge-Kutta method of order $M+1$ \cite{Gottlieb2001} and the Rusanov (local Lax-Friedrichs) numerical flux is used.

\def\arraystretch{1.1}
\begin{table}
{\scriptsize
\begin{tabular}{r r | c r | c r | c r | c r |}
 & &  \multicolumn{4}{|c|}{$N=1$} & \multicolumn{4}{|c|}{$N=2$} \\ 
 & &   \multicolumn{2}{|c}{Limiter \emph{off}} & \multicolumn{2}{c|}{Limiter \emph{on}} & \multicolumn{2}{|c}{Limiter \emph{off}} & \multicolumn{2}{c|}{Limiter \emph{on}} \\ \hline
 \multirow{7}{*}{$M=2$}	&$I$ & $E_2$ &  $\mathcal{O}_2 $	& $E_2$ &  $\mathcal{O}_2 $ &  $E_2$ & $\mathcal{O}_2 $	& $E_2$ &  $\mathcal{O}_2 $
	 \\ \hline 
&	$10$  &  $2.42 \text{E}-01$ &  		   & $2.54\text{E}-01$ &		   & & & & \\
&	$20$  &  $3.80 \text{E}-02$ & $2.67$   & $6.36\text{E}-02$ & $1.99$  & & & & \\
&	$40$  &  $3.95 \text{E}-03$ & $3.27$   & $1.43\text{E}-02$ & $2.15$  & & & & \\
&	$80$  &  $5.26 \text{E}-04$ & $2.91$   & $5.24\text{E}-03$ & $1.45$  & & & & \\
&	$160$ &  $8.43 \text{E}-05$ & $2.64$   & $1.03\text{E}-03$ & $2.35$  & & & & \\
	\hline 
 \multirow{5}{*}{$M=3$}	
&	$10$  &  $2.03 \text{E}-01$ &  		    & $2.32\text{E}-01$ &  			&  $1.66 \text{E}-02$ &  		& $2.61 \text{E}-02$ & 		\\
&	$20$  &  $2.23 \text{E}-02$ & $3.19$  	& $5.26\text{E}-02$ & $2.14$ 	&  $6.00 \text{E}-04$ & $4.79$	& $4.04 \text{E}-03$ & $2.42$\\
&	$40$  &  $1.25 \text{E}-03$ & $4.16$  	& $9.39\text{E}-03$ & $2.49$	&  $3.63 \text{E}-05$ & $4.04$ 	& $4.87 \text{E}-04$ & $3.05$\\ 
&	$80$  &  $7.56 \text{E}-05$ & $4.05$  	& $2.99\text{E}-03$ & $1.65$	&  $2.85 \text{E}-06$ & $3.67$  & $3.70 \text{E}-05$ & $3.72$\\
&	$160$ &  $5.63 \text{E}-06$ & $3.75$  	& $9.91\text{E}-04$ & $1.59$	&  $2.43 \text{E}-07$ & $3.55$  & $2.43 \text{E}-07$ & $7.25$\\
	\hline
 \multirow{5}{*}{$M=4$}	
&	$10$  &  $6.93 \text{E}-02$ &  			& $1.91 \text{E}-01$ &  		&  $1.59 \text{E}-02$  &  		 & $2.08 \text{E}-02$ &  	  \\
&	$20$  &  $1.95 \text{E}-03$ & $5.15$  	& $5.35 \text{E}-02$ & $1.84$ 	&  $5.04 \text{E}-04$  & $4.98$  & $3.73 \text{E}-03$ & $2.48$\\
&	$40$  &  $7.87 \text{E}-05$ & $4.63$  	& $1.42 \text{E}-02$ & $1.91$	&  $2.11 \text{E}-05$  & $4.58$  & $4.31 \text{E}-04$ & $3.12$\\ 
&	$80$  &  $5.76 \text{E}-06$ & $3.78$  	& $3.39 \text{E}-03$ & $2.06$	&  $9.36 \text{E}-07$  & $4.49$  & $3.14 \text{E}-05$ & $3.78$\\
&	$160$ &  $4.95 \text{E}-07$ & $3.54$  	& $6.86 \text{E}-04$ & $2.31$	&  $4.15 \text{E}-08$  & $4.49$  & $4.15 \text{E}-08$ & $9.57$\\
	\hline
 \multirow{5}{*}{$M=5$}	
&	$10$  &  $4.20 \text{E}-02$ &  			& $1.78 \text{E}-01$ &  		&  $3.34 \text{E}-03$  &    	 & $1.61 \text{E}-02$ &  	  \\
&	$20$  &  $4.50 \text{E}-04$ & $6.55$  	& $5.24 \text{E}-02$ & $1.76$  	&  $3.72 \text{E}-05$  & $6.49$  & $4.51 \text{E}-03$ & $1.83$\\
&	$40$  &  $4.65 \text{E}-05$ & $3.27$  	& $1.51 \text{E}-02$ & $1.80$	&  $7.04 \text{E}-07$  & $5.72$  & $4.91 \text{E}-04$ & $3.20$\\ 
&	$80$  &  $5.24 \text{E}-06$ & $3.15$  	& $4.07 \text{E}-03$ & $1.87$	&  $1.52 \text{E}-08$  & $5.54$  & $3.65 \text{E}-05$ & $3.75$\\
&	$160$ &  $4.83 \text{E}-07$ & $3.44$  	& $9.09 \text{E}-04$ & $2.16$	&  $3.33 \text{E}-10$  & $5.51$  & $3.33 \text{E}-10$ & $16.74$\\
	\hline
 \multirow{5}{*}{$M=6$}	
&	$10$  &  & 	&  & 														&  $2.68 \text{E}-03$  &    	& $1.62 \text{E}-02$ & 		 \\
&	$20$  &  & 	&  &  														&  $2.00 \text{E}-05$  & $7.07$ & $4.39 \text{E}-03$ & $2.54$\\
&	$40$  &  & 	&  &  														&  $2.66 \text{E}-07$  & $6.23$ & $4.72 \text{E}-04$ & $3.71$\\ 
&	$80$  &  & 	&  &  														&  $4.47 \text{E}-09$  & $5.89$ & $3.43 \text{E}-05$ & $3.89$\\
&	$160$ &  & 	&  &  														&  $9.00 \text{E}-11$  & $5.64$ & $9.00 \text{E}-11$ & $18.54$\\
 \hline
\end{tabular}
}
\end{table}

\begin{table}
{\scriptsize
\begin{tabular}{r  r | c r | c r | c r | c r |}
 & &  \multicolumn{4}{|c|}{$N=3$} & \multicolumn{4}{|c|}{$N=4$} \\ 
 & &   \multicolumn{2}{|c}{Limiter \emph{off}} & \multicolumn{2}{c|}{Limiter \emph{on}} & \multicolumn{2}{|c}{Limiter \emph{off}} & \multicolumn{2}{c|}{Limiter \emph{on}} \\ \hline
 \multirow{7}{*}{$M=4$}	&$I$ & $E_2$ &  $\mathcal{O}_2 $	& $E_2$ &  $\mathcal{O}_2 $ &  $E_2$ & $\mathcal{O}_2 $	& $E_2$ &  $\mathcal{O}_2 $
	 \\ \hline 
&	$10$  &  $1.07 \text{E}-03$ &  			& $1.36 \text{E}-03$ & 			& & & & \\
&	$20$  &  $3.06 \text{E}-05$ & $5.13$  	& $1.78 \text{E}-04$ & $2.93$ 	& & & & \\
&	$40$  &  $1.04 \text{E}-06$ & $4.89$  	& $1.04 \text{E}-06$ & $7.43$	& & & & \\
&	$80$  &  $4.08 \text{E}-08$ & $4.67$  	& $4.08 \text{E}-08$ & $4.67$	& & & & \\
&	$160$ &  $1.74 \text{E}-09$ & $4.55$  	& $1.74 \text{E}-09$ & $4.55$	& & & & \\
	\hline
 \multirow{5}{*}{$M=5$}	
&	$10$  &  $1.00 \text{E}-03$ &  			& $1.38 \text{E}-03$ &  		&  $9.79 \text{E}-05$ &  		    & $1.65 \text{E}-04$ 	&  	\\
&	$20$  &  $2.42 \text{E}-05$ & $5.38$  	& $1.79 \text{E}-04$ & $2.95$  	&  $1.48 \text{E}-06$ & $6.05$   	& $1.48 \text{E}-06$ 	& $6.80$\\
&	$40$  &  $5.59 \text{E}-06$ & $5.43$  	& $5.59 \text{E}-07$ & $8.32$	&  $2.55 \text{E}-08$ & $5.86$  	& $2.55 \text{E}-08$ 	& $5.86$\\ 
&	$80$  &  $1.25 \text{E}-08$ & $5.49$  	& $1.25 \text{E}-08$ & $5.49$	&  $5.14 \text{E}-10$ & $5.63$  	& $5.14 \text{E}-10$ 	& $5.63$\\
&	$160$ &  $2.76 \text{E}-10$ & $5.50$  	& $2.76 \text{E}-10$ & $5.50$	&  $1.11 \text{E}-11$ & $5.53$  	& $1.11 \text{E}-11$ 	& $5.53$\\
	\hline
 \multirow{5}{*}{$M=6$}	
&	$10$  &  $1.84 \text{E}-04$ &  			& $1.73 \text{E}-03$ &  		&  $8.94 \text{E}-05$ &  		    & $ 1.88 \text{E}-04$ 	&  	\\
&	$20$  &  $2.71 \text{E}-06$ & $6.09$  	& $1.97 \text{E}-04$ & $3.13$  	&  $1.11 \text{E}-06$ & $6.33$   	& $ 1.11 \text{E}-06$ 	& $7.40$\\
&	$40$  &  $2.80 \text{E}-08$ & $6.59$  	& $2.80 \text{E}-08$ & $12.78$	&  $1.28 \text{E}-08$ & $6.45$   	& $ 1.28 \text{E}-08$ 	& $6.45$\\ 
&	$80$  &  $1.59 \text{E}-10$ & $7.46$  	& $1.59 \text{E}-10$ & $7.46$	&  $1.43 \text{E}-10$ & $6.48$  	& $ 1.43 \text{E}-10$ 	& $6.48$\\
&	$160$ &  $3.14 \text{E}-12$ & $5.66$  	& $3.14 \text{E}-12$ & $5.66$	&  $1.64 \text{E}-12$ & $6.45$  	& $ 1.64 \text{E}-12$ 	& $6.45$\\
	\hline
\end{tabular}
}
\caption{$L^2$ errors for the linear advection equation with smooth initial data}
\end{table}
We observe that on coarse grids the limited version of the $P_NP_M$ scheme is only of order $N+1$. On finer grids, however, the limiter is never active and the full order $M+1$ of the unlimited scheme is achieved. At which level of grid refinement the limiter becomes inactive depends on $N$: For larger $N$, no limiting is needed on rather coarse grids, e.g. for $P_4P_5$ and $P_4P_6$ the limiter is inactive on a grid with $20$ cells, which corresponds to a grid size of $h=1/10$. Generally speaking, if the representation of the data in the lower order polynomials is already sufficiently good, then the high order representation does not need to be limited.
\subsection{Burgers equation}
Consider Burgers equation
\begin{equation}
\label{Burgers CL}
\frac{\partial u}{\partial t} + \frac{\partial }{\partial x} \left( \frac{u^2}{2}\right) = 0,
\end{equation}
with initial data
\begin{equation}
\label{Burgers IC}
	u(x,0) = - 5 \exp \left( -50\left( x-\frac{1}{2} \right)^2 \right) 
			 + 5 \exp \left( -50\left( x+\frac{1}{2} \right)^2 \right). 
\end{equation}
We solve \eqref{Burgers CL}, \eqref{Burgers IC} on $(-1,1)$ with transmissive boundary conditions. We use a fourth order SSP Runge-Kutta method and denote by $\theta$ the average value of the limiter during the Runge-Kutta stages for one time-step. We present the results for $P_2P_4$ on $160$ cells and for $P_4P_6$ at several time points during the simulation.

\begin{figure}
\resizebox{\hsize}{!}{		
	\includegraphics*[height=0.2\textheight]{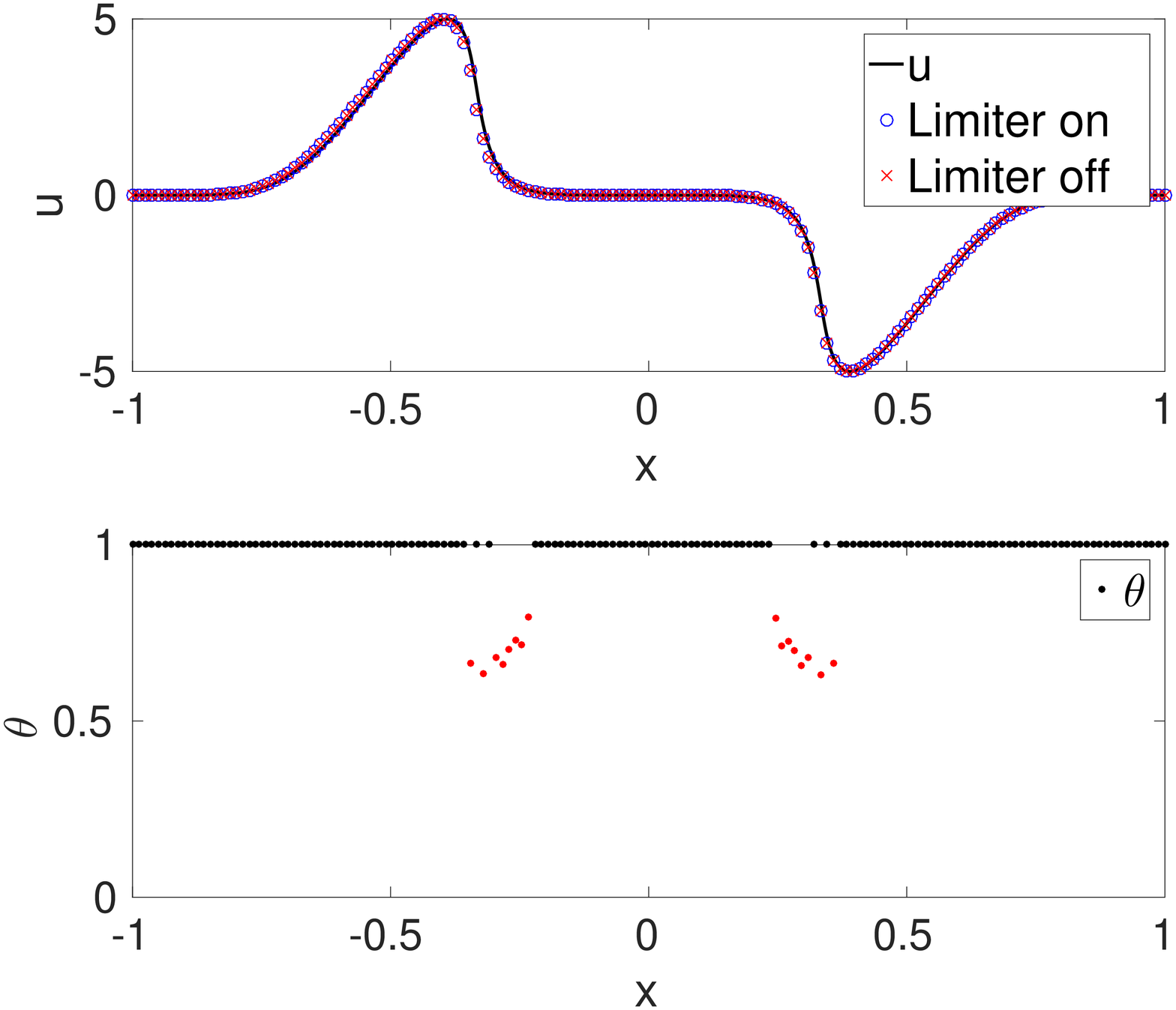}
	\includegraphics*[height=0.2\textheight]{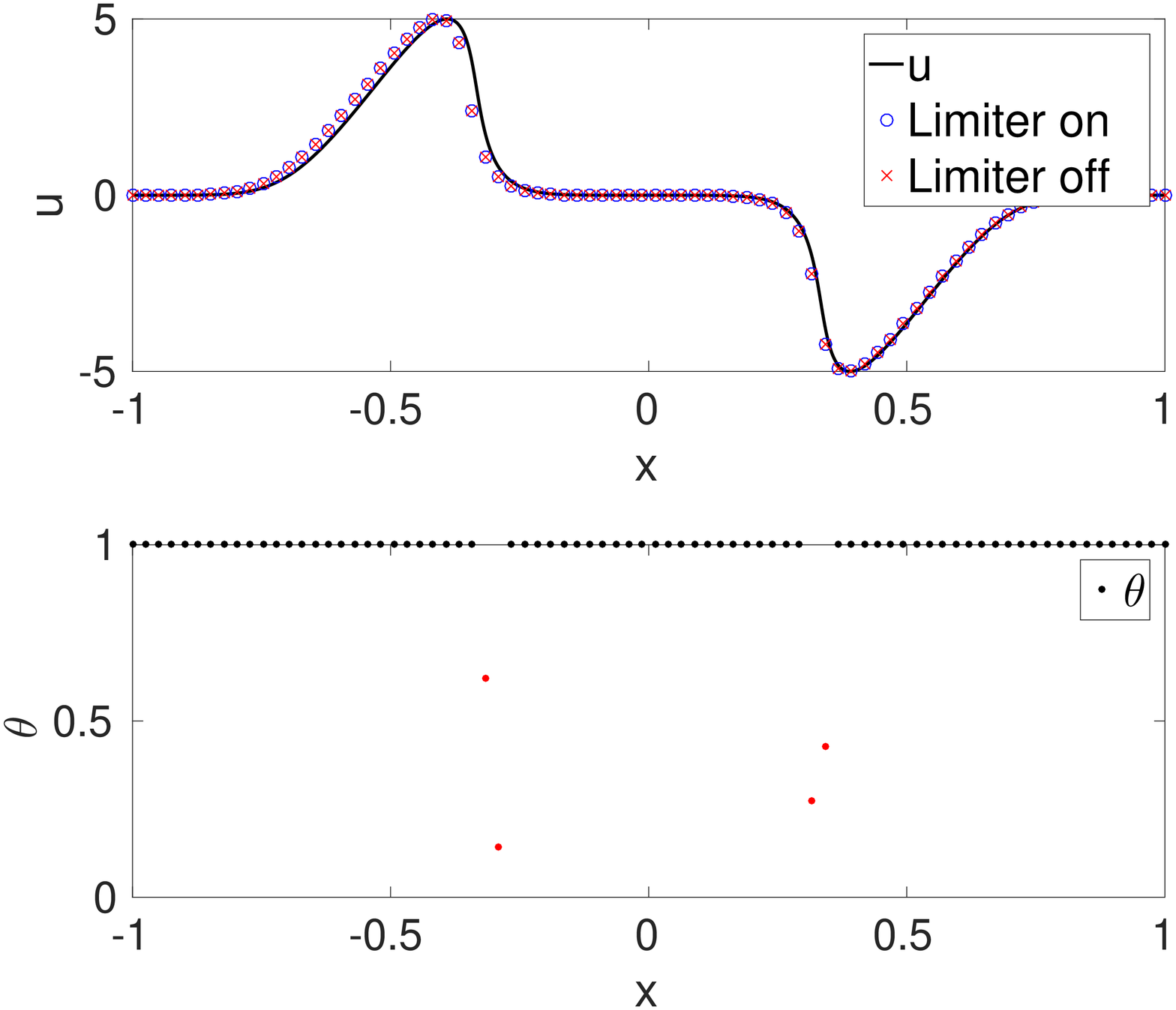}
}	

\smallskip

\resizebox{\hsize}{!}{	
	\includegraphics*[height=0.2\textheight]{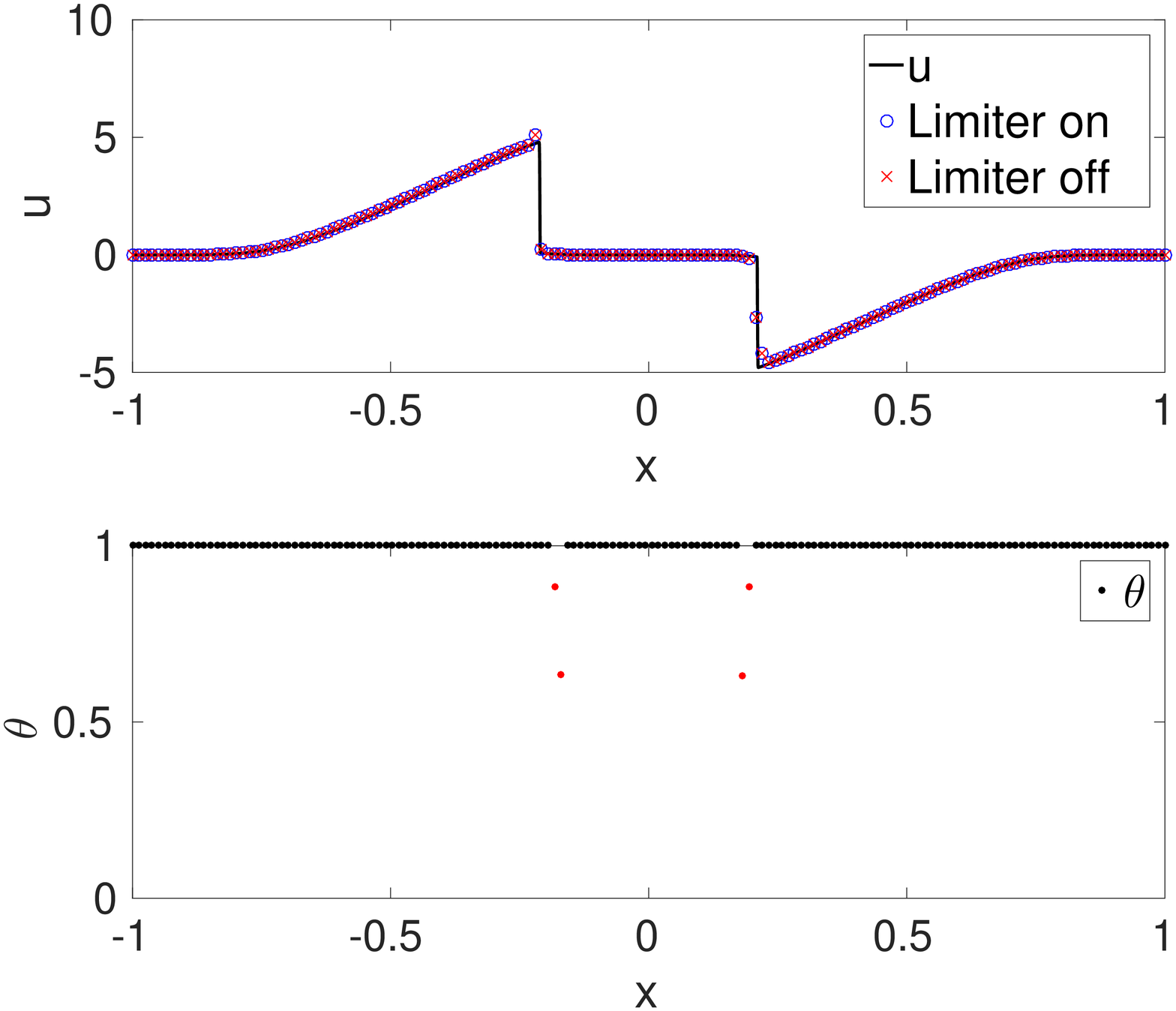}
	\includegraphics*[height=0.2\textheight]{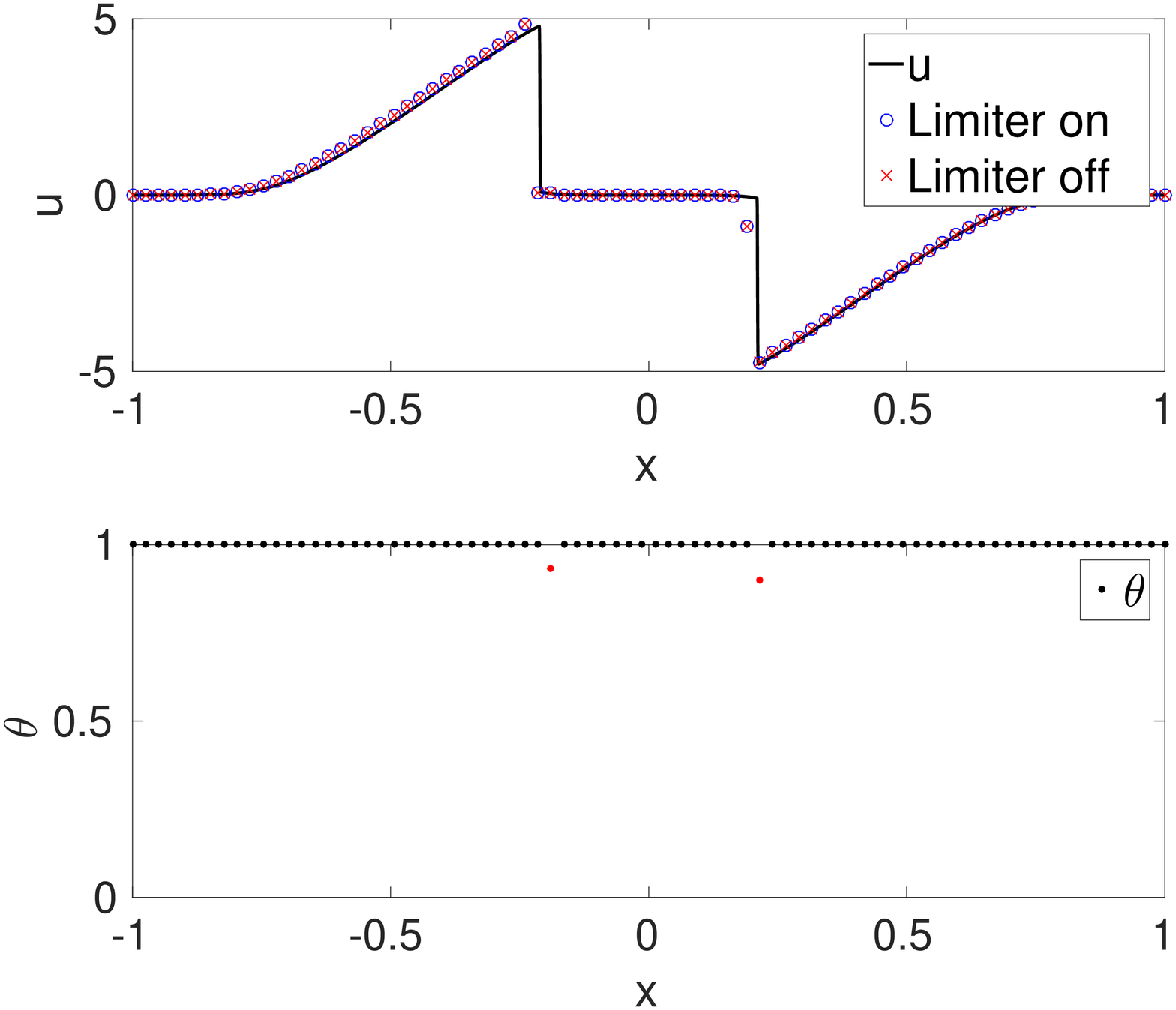}
}
\smallskip

\resizebox{\hsize}{!}{	
	\includegraphics*[height=0.2\textheight]{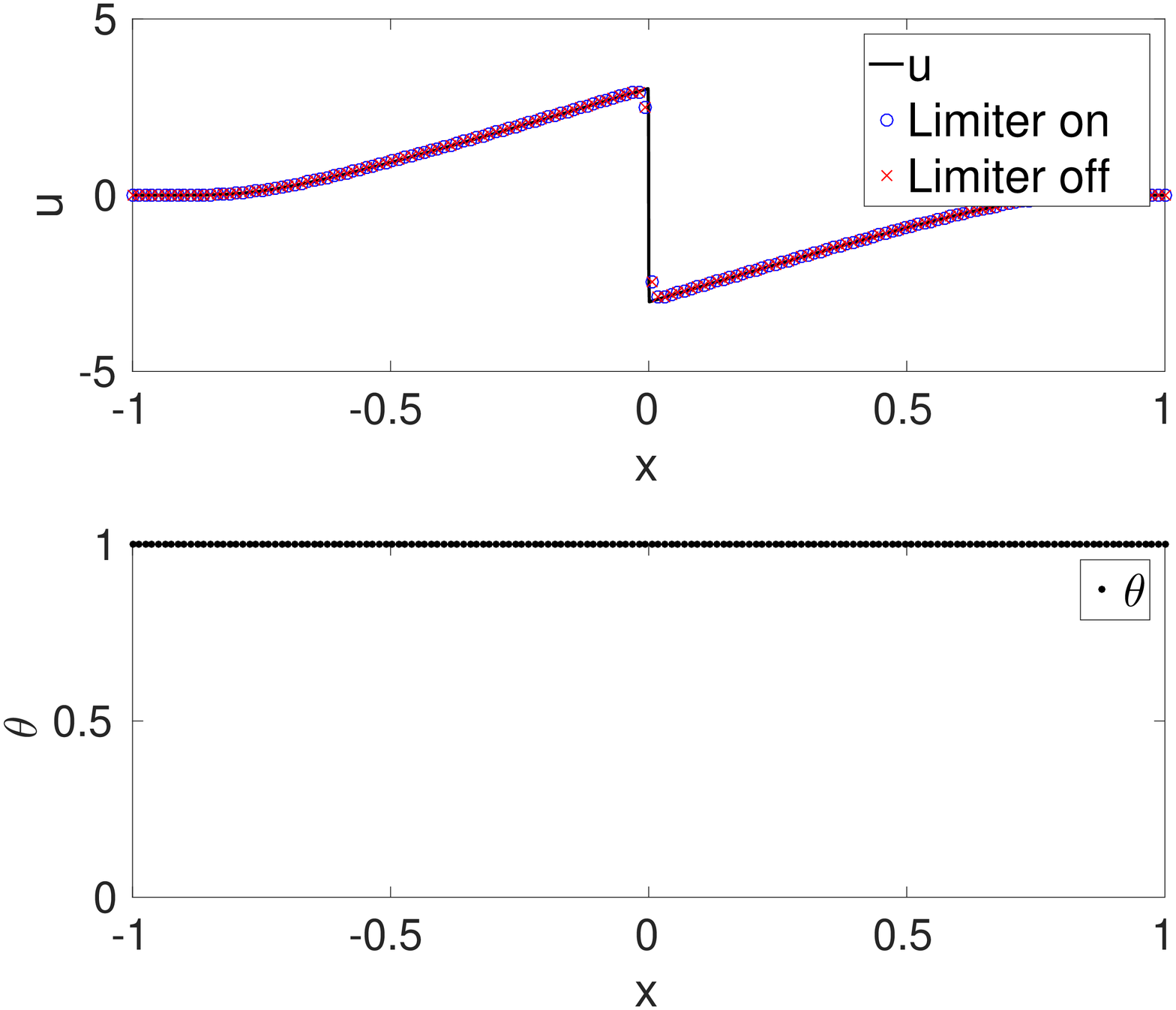}
	\includegraphics*[height=0.2\textheight]{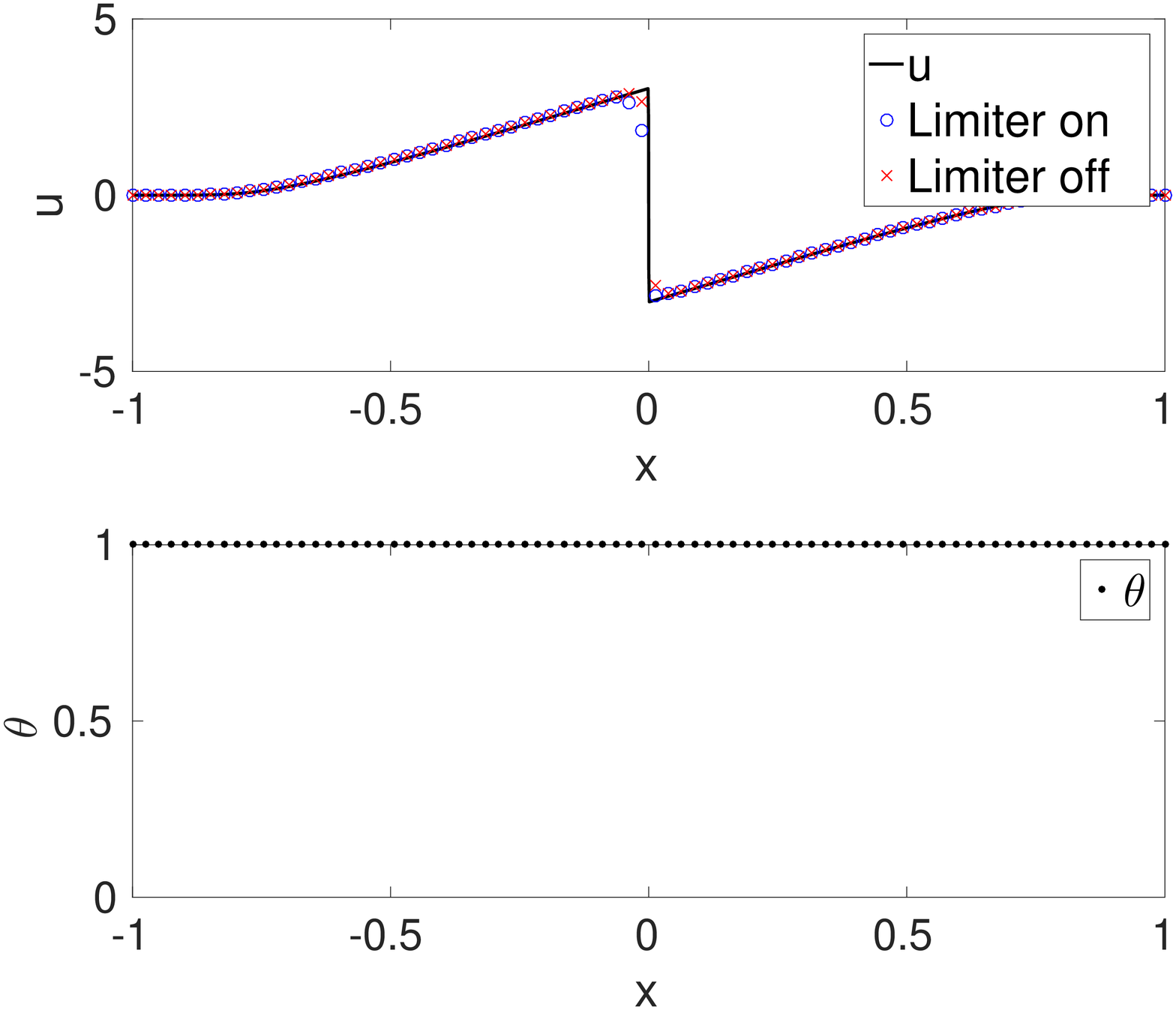} 
}
\begin{center}
 c) $ t = 0.6$
\end{center}

\caption{Burgers equation, left: $P_2P_4$ solution on $160$ cells, right: $P_4P_6$ solution on $80$ cells. Top row: $t=0.022$, middle row: $t=0.066$, bottom row: $t=0.198$ \label{fig:burgers} }
\end{figure}

Results are shown in Figure \ref{fig:burgers}. We observe a similar behaviour as in the linear case: In regions where the solution is smooth, the limiter is not active. In fact, most limiting is needed shortly before the shocks form. Again we can see that if the $P_N$ representation of the solution is sufficiently accurate, even on rather coarse grid the $P_M$ part requires no or only very little limiting.

\subsection{Traffic flow}
Consider the following Lighthill-Whitman type model for traffic flow:
\begin{equation}
\label{traffic CL}
	\frac{\partial \rho}{\partial t} + \frac{\partial}{\partial x}\left( 2\rho \exp \left(-\frac{1}{2}\rho^2 \right)\right)  = 0,
\end{equation}
where $\rho \in (0,1)$ describes the density of cars on a road. Note that in this case $f$ is a \emph{strictly concave} flux.  We solve \eqref{traffic CL} on $(-1,1)$ with periodic boundary conditions and initial data given by
\[
	\rho(x,0) = \frac{1}{2} + \frac{1}{4}\sin (\pi x).
\]
The solution is given by a sinusodial wave travelling to the right that is deformed until a shock emerges. Note that in this model, lower values of $\rho$ lead to a faster speed of propagation.

\begin{figure}
\resizebox{\hsize}{!}{		
	\includegraphics*[height=0.2\textheight]{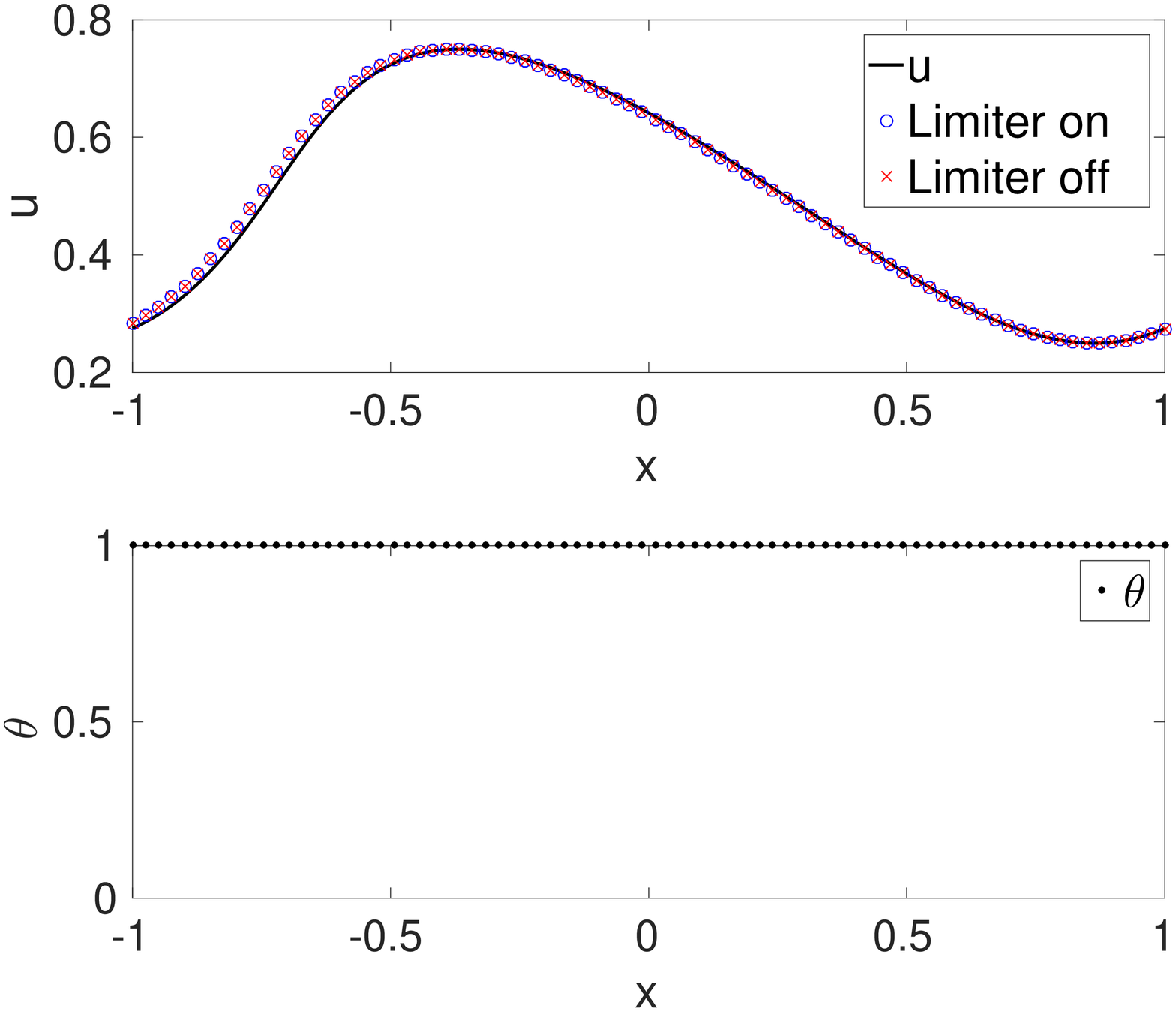}
	\includegraphics*[height=0.2\textheight]{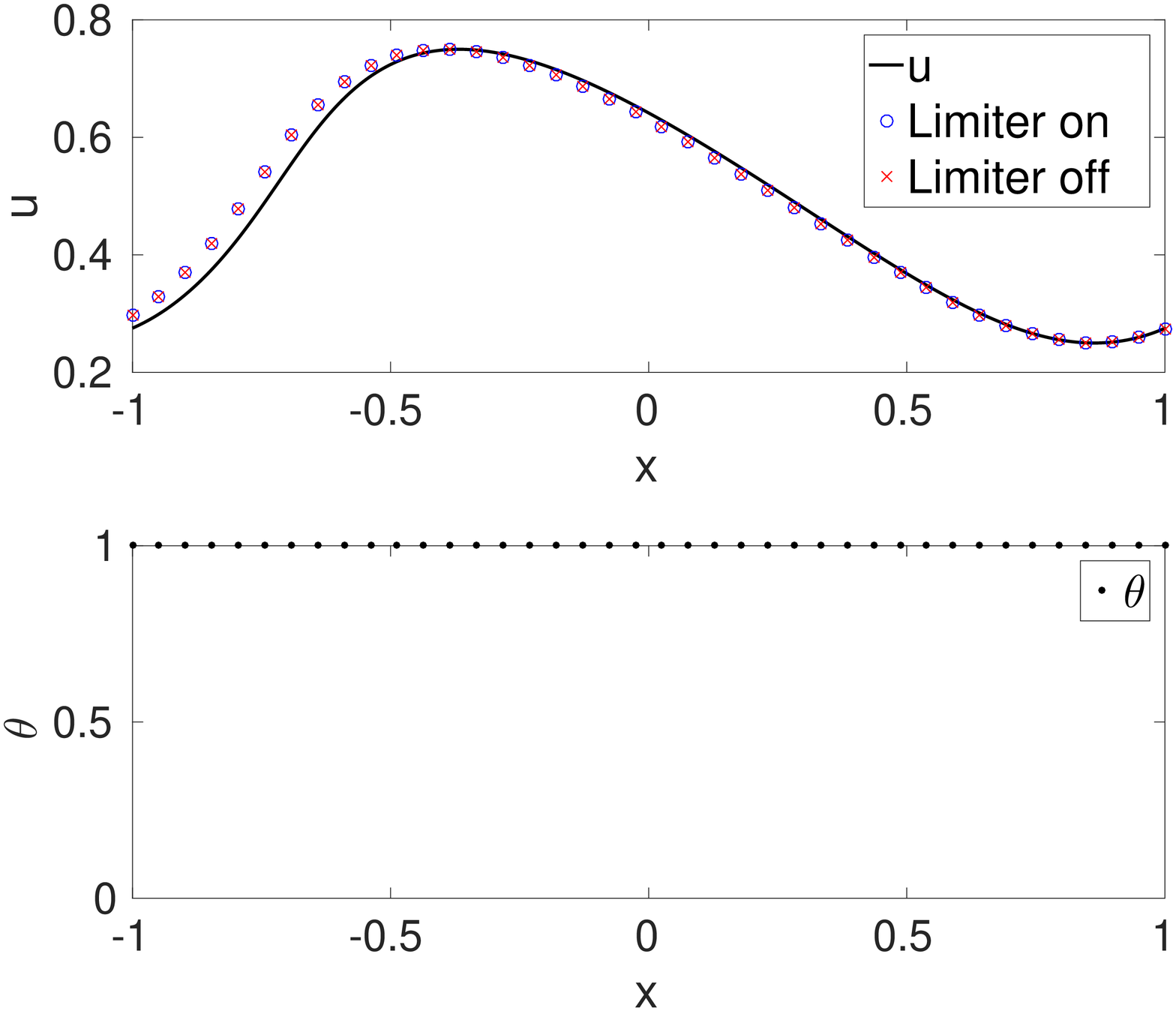}
}	

\smallskip

\resizebox{\hsize}{!}{	
	\includegraphics*[height=0.2\textheight]{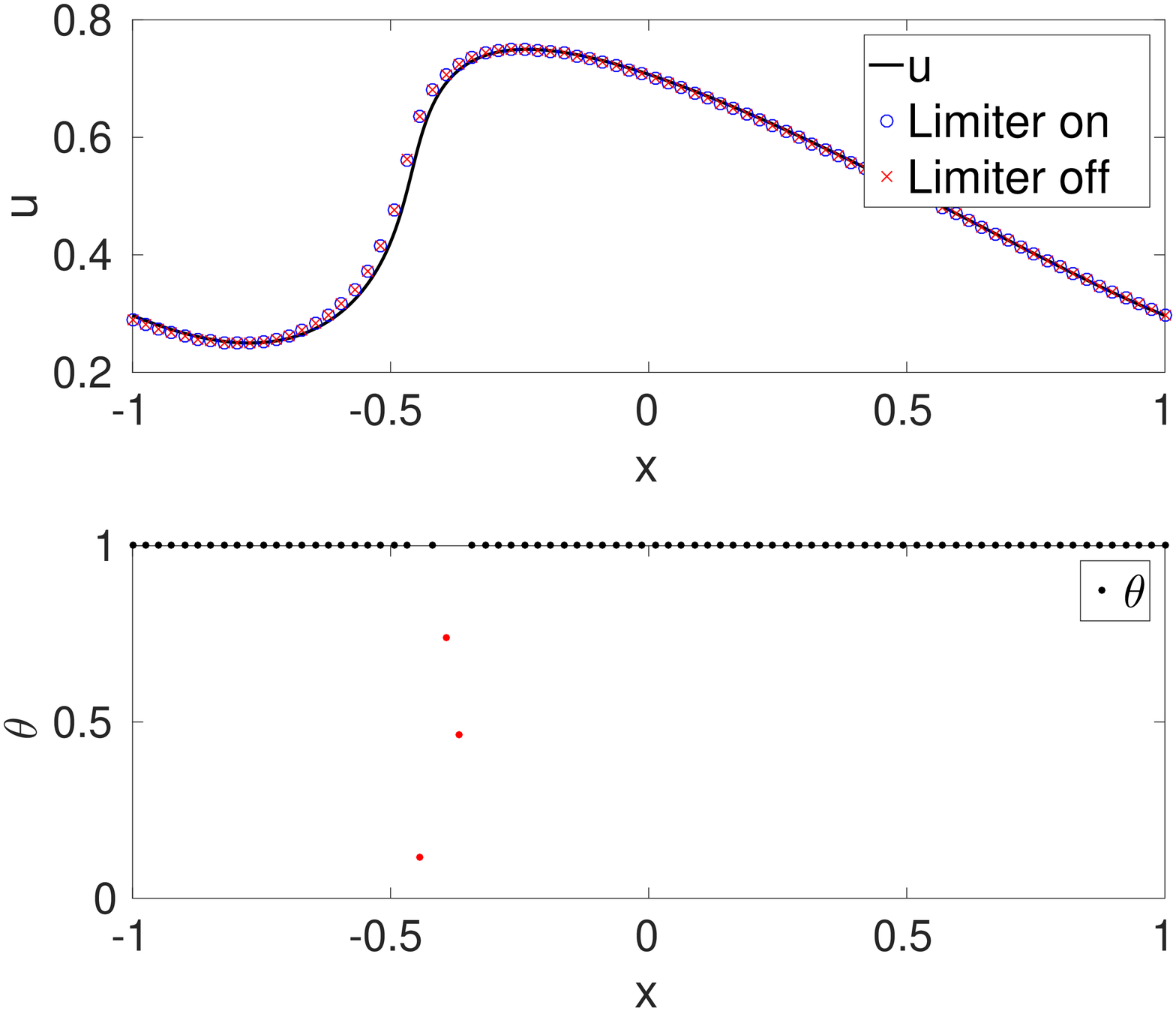}
	\includegraphics*[height=0.2\textheight]{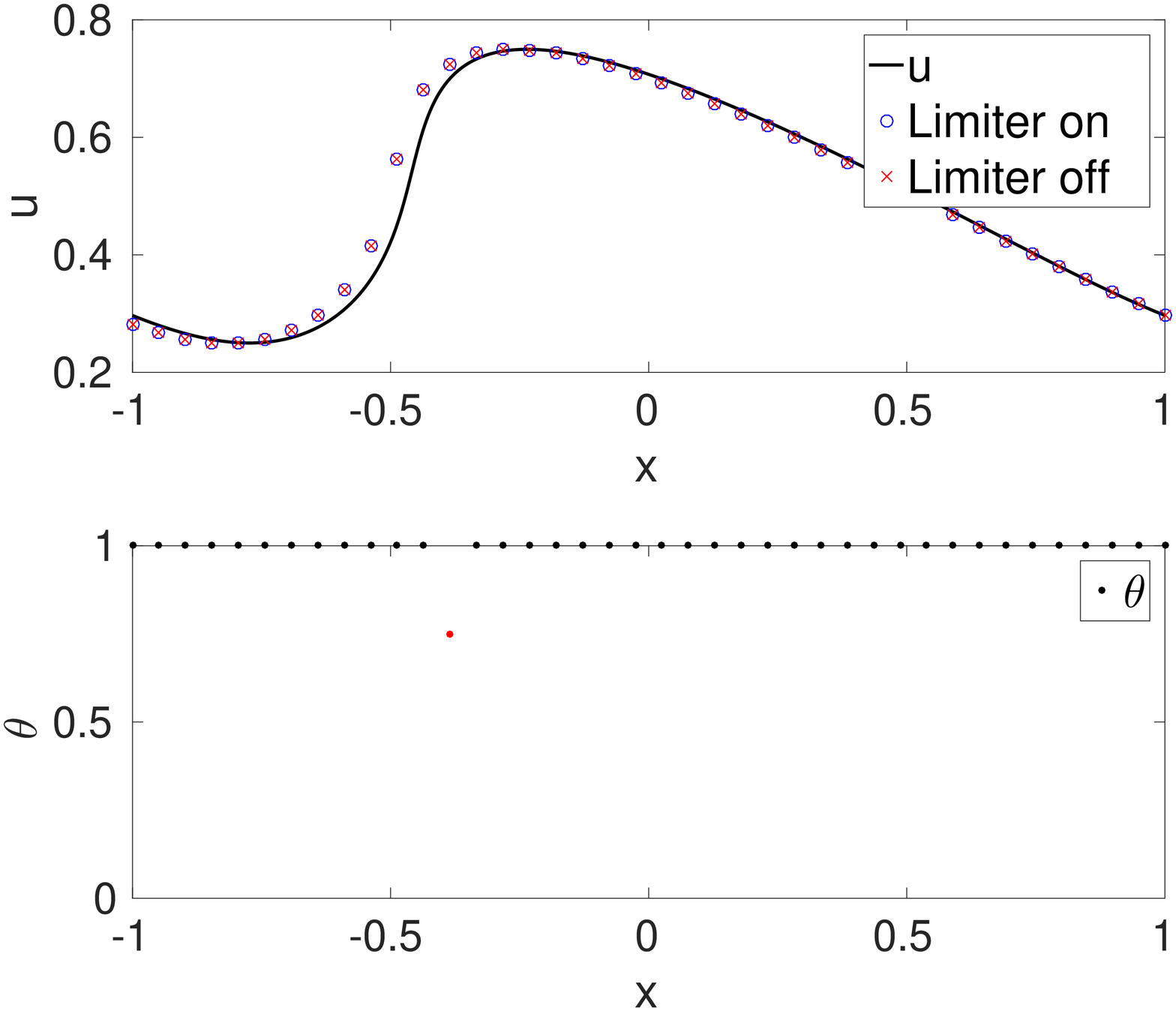}
}
\smallskip

\resizebox{\hsize}{!}{	
	\includegraphics*[height=0.2\textheight]{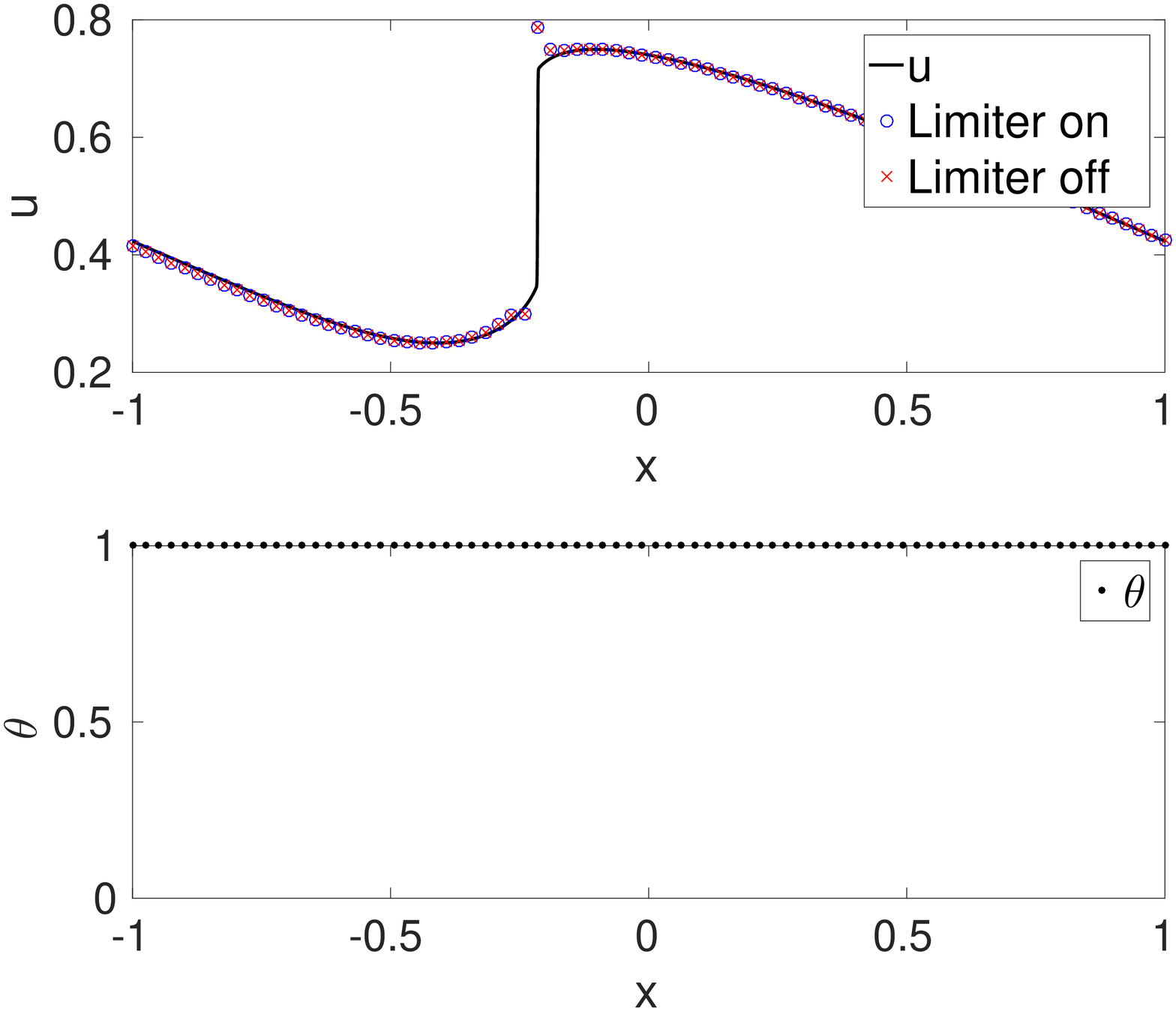}
	\includegraphics*[height=0.2\textheight]{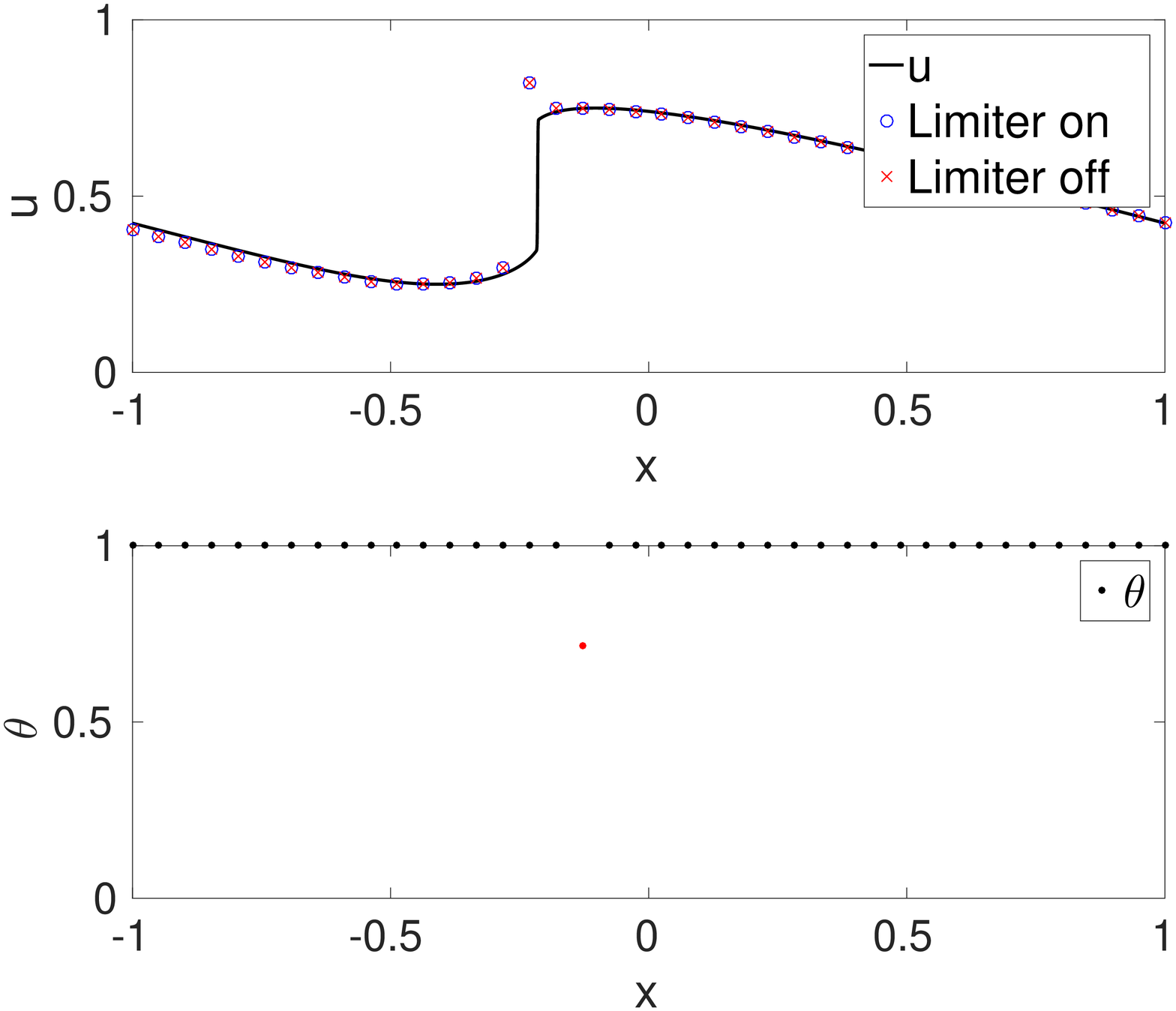} 
}
\begin{center}
 c) $ t = 0.6$
\end{center}

\caption{Traffic flow, left: $P_2P_4$ solution on $80$ cells, right: $P_4P_6$ solution on $40$ cells. Top row: $t=0.2$, middle row: $t=0.4$, bottom row: $t=0.6$ \label{fig:traffic} }
\end{figure}

Numerical results are shown in Figure \ref{fig:traffic}. The left column depicts the numerical solution on $80$ cells and the left column shows the results for $P_4P_6$ on $40$ cells. Again, time integration is performed with a fourth order SSP Runge-Kutta method and $\theta$ denotes the average value of the limiter during the respective time-step. We observe that the limiter is only active on a very small number of cells, shortly before the shock is formed and close to the shock front after its formation, while the numerical solution is not limited as long as the exact solution is smooth. 

\section{Conclusion}\label{sec:conclusion}

We have seen that the reconstruction step in $P_NP_M$ schemes can cause unstable behaviour with respect to the $L^2$ norm of the numerical solution. We demonstrated analytically how a numerical flux using reconstructed function values can lead to a violation of a cell entropy condition. In particular for linear problems we were able to derive a simple algebraic relation between the jump in the reconstruction and the jump in the data that determines whether the scheme is square entropy stable, or not. 

With this new condition, that also has a nonlinear analogon, we were able to construct a flux limiter that guarantees the entropy stability of the numerical solution. As expected, entropy stability only becomes an issue when the solution contains large jumps or very steep gradients, in which case the reconstruction can become oscillatory. Numerical experiments validate the limiter's ability to stabilize the method while still maintaining high order accuracy in smooth regions for sufficiently fine meshes.

\appendix
\section{Existence and uniqueness of the reconstruction in the case $M=3N+2$.}

We show the existence and uniqueness of the solution to the reconstruction problem
\begin{equation}
	\sum_{m=0}^{3N+2} \hat{w}_m \langle \Psi_k, \Phi_\ell^{(j)} \rangle_{T^{(j)}} = \sum_{m=0}^{N} \hat{v}_k^{(j)} \langle  \Phi_k^{(j)},  \Phi_\ell^{(j)} \rangle_{T^{(j)}}, 
\end{equation}
for $ \ell=0,\dots, N, ~ j\in \{i-1,i,i+1\}$. Recall that the functions $\{\Phi_\ell^{(j)}~:~ \ell = 0,\dots,N\}$ are the shifted Legendre polynomials of degree $\ell$ on $T^{(j)}$ and the $\{ \Psi_m~:~m=0,\dots 3N+2 \}$ are Legendre polynomials on the central cell $T^{(i)}$, extended to the whole stencil.

A basis $\{ \Theta_k~:~k=0,\dots, 3N+2\}$ of $V_h^N\left( \mathcal{S}^{(i)} \right)$, that is orthogonal in the sense that
\[
	\left\langle \Theta_m, \Theta_n \right\rangle_{\mathcal{S}^{(i)}} = 
	\int_{\mathcal{S}^{(i)} } \Theta_m(x) \Theta_n(x)~dx = 0, 
	\qquad \text{ if } m \neq n.
\]
is given by
\[
	\tilde{\Phi}_\ell^{(j)}(x) =~ \left\{
	\begin{array}{cc}
		\Phi_\ell^{(j)}(x), & \quad x \in T^{(j)},\\
		0, 					& \quad \text{ else }.
	\end{array}
	\right.
\]
and an an index $k \equiv k(\ell, j)$, such that
\begin{align*}
	k(0,i-1) & = 0,    \quad	\dots, 	& k(N, i-1)	& = N,  &\\ 
	k(0, i ) & = N+1,  \quad 	\dots, 	& k(N, i )	& = 2N+1, &\\
	k(0,i+1) & = 2N+2, \quad	\dots, 	& k(N, i+1) & = 3N+2. &
\end{align*}
We let 
\[
	\Theta_k = \tilde{\Phi}_\ell^{(j)}, \qquad \ell \in \{0, \dots, N\}, ~ j \in \{i-1, i, i+1\},
	\qquad  k = k(\ell,j).
\]
With this, the reconstruction problem can be written as
\begin{equation}\label{inv L2 appendix}
	\sum_{k=0}^{3N+2} \hat{w}_k \langle \Psi_k, \Theta_m \rangle_{\mathcal{S}^{(i)}} = \sum_{k=0}^{3N+2} \hat{v}_k \langle  \Theta_k,  \Theta_m \rangle_{\mathcal{S}{(i)}}, \quad m = 0,\dots, 3N+2,
\end{equation}
where we denote $\hat{v}_k \equiv \hat{v}_{k(\ell,j)} = \hat{v}_\ell^{(j)}$. It is worth pointing out that condition \eqref{inv L2 appendix} implies that $v$ is the $L^2$-projection of $w$ onto $V_h^N\left( \mathcal{S}^{(i)} \right)$. Thus, our reconstruction problem is an \emph{inverse projection} problem, for which existence and uniqueness of the solution are non-trivial. 

In order to guarantee existence and uniqueness for the solution of the reconstruction problem, we have to check whether the matrix
\[		
	\tilde{A} =  \left( \langle\Theta_m, \Psi_\ell \rangle_{\S} \right)_{m,\ell \in \{0, \dots, 3N+2\}} 
\]
is invertible. 

At first, note that the conditions for the reconstruction are invariant under linear coordinate-transformations. So we can map $T^{(i-1)}$ to the interval $(-3,-1)$, $T^{(i)}$ to $(-1,1)$, and $T^{(i+1)}$ to $(1,3)$. Denote the Legendre polynomial of degree $k$ on $(-1,1)$ by $P_k$ and the the shifted Legendre polynomials of degree $k$ on $(-3,1)$ and $(1,3)$ by $P_k^{(-1)}$ and $P_k^{(+1)}$, respectively. 

Then, after the coordinate change, the basis functions $\Phi_k^{(i-1)}$ and $\Phi_k^{(i+1)}$ become $P_k^{(-1)}$ and $P_k^{(+1)}$, respectively, while the 
$ \Phi_k^{(i)}$ simply become $P_k$. Moreover, the $\Psi_k$ also become $P_k$. Thus, on $(-3,-1)$ we get the matrix coefficient
\begin{align*}
	\int_{-3}^{-1}  P_k(s) P_\ell^{(-1)}(s)~ds  
& 	= \int_{-1}^{1} P_k(s+2) P_\ell^{(-1)}(s+2)~ds\\  
&	= \int_{-1}^{1} P_k(s+2) P_\ell(s)~ds,
\end{align*}
by shifting to $(-1,1)$ and noting that by construction $P_\ell^{(-1)}(s+2) = P_\ell(s)$. Similarly, we have
\[
	\int_{1}^{3}  P_\ell(s) P_k^{(+1)}(s)~ds = \int_{-1}^{1} P_\ell(s-2) P_k(s)~ds.
\]
By a simple symmetry argument, it is then easy to see that
\[
	\int_{-1}^{1} P_\ell(s+2) P_k(s)~ds = (-1)^{k+\ell} \int_{-1}^{1} P_\ell(s-2) P_k(s)~ds.
\]

Let 
\[
	a_{k,\ell} = \int_{-1}^1 P_k(s)P_\ell(s+2)~ds, \qquad 	b_{k,\ell} = (-1)^{k+\ell}a_{k, \ell}
\]
Then the matrix $\tilde{A}$ is invertible, if and only if the matrix
\[
A =	\left[
		\arraycolsep=3pt\def\arraystretch{2.5}
		\begin{array}{c}
			\displaystyle 	\int_{-1}^1 P_k(s) P_\ell(s)~ds \medskip \\ 
			\hdashline
			\displaystyle 	\int_{-1}^1 P_k(s) P_\ell(s+2)~ds \medskip \\
			\hdashline
			\displaystyle 	\int_{-1}^1 P_k(s) P_\ell(s-2)~ds \\
		\end{array}
	\right]_{\scriptsize\begin{aligned} k &= 0, \dots, N \\[-4pt] \ell & = 0, \dots, 3N+2\end{aligned}}
=		
\left[\arraycolsep=3pt\def\arraystretch{2.5}
		\begin{array}{c}
			\displaystyle 	\frac{1}{2k+1}\delta_{k,\ell} \medskip \\ 
			\hdashline
			\displaystyle 	a_{k,\ell} \medskip \\
			\hdashline
			\displaystyle 	b_{k, \ell} \\
		\end{array}
	\right]_{\scriptsize\begin{aligned} k &= 0, \dots, N \\[-4pt] \ell & = 0, \dots, 3N+2\end{aligned}} 
\]
is invertible. We have 
\[
	A = 
	\left[
		\begin{array}{cc}
		D & 0\\
		A_{21} & A_{22}
		\end{array}
	\right]
\]
with a diagonal part $D \in \mathbb{R}^{(N+1)\times(N+1)}$ and
\begin{align*}
A_{21} & = \left[
	\def\arraystretch{2.5}
	\begin{array}{rcr}
		\displaystyle a_{k,0}  & \dots & 	a_{k,N}  \medskip \\
		\hdashline
		\displaystyle (-1)^k a_{k,0} 	& \dots 					& (-1)^{k+N} a_{k,N}
	\end{array}
\right]_{k=0,\dots,N} \in \mathbb{R}^{(2N+2)\times(N+1)}\\
A_{22} & = \left[
	\def\arraystretch{2.5}	
	\begin{array}{rcr}
		 a_{k,N+1} 				& \dots  &	a_{k,3N+2} \medskip \\
		 \hdashline
		 (-1)^{k+N+1} a_{k,N+1} & \dots  & (-1)^{k+3N+2}a_{k,3N+2}
	\end{array}
	\right]_{k=0,\dots,N} \in \mathbb{R}^{(2N+2)\times(2N+2)}
\end{align*}
Defining 
\[
	c_{k,\ell} = \frac{1}{2}\left( a_{k,\ell} + b_{k,\ell} \right),
\]
it is then sufficient to show that
\[
	 B = \left[
		\arraycolsep=1pt\def\arraystretch{1}
		\begin{array}{ccc}
			a_{0,N+1} 	& \dots & a_{0, 3N+3} \\ 
			\vdots 		&~		& \vdots \\ 
			a_{N,N+1}	& \dots & a_{N, 3N+3} \medskip \\ 
			\hdashline
			c_{0,N+1} 	& \dots & c_{0, 3N+3} \\ 
			\vdots 		&~		& \vdots \\ 
			c_{N,N+1}	& \dots & c_{N, 3N+3}
		\end{array}
	\right] \in \mathbb{R}^{(2N+2)\times (2N+2)}
\]
is invertible. Before we show this, let us proof the following:

{\bf Lemma.} {\it Let $0 \neq Q $ be polynomial of degree $N$ that does not vanish identically on $(-1,1)$, such that for some $\ell > N$ we have
\[
	\int_{-1}^1 P_\ell(s+2)Q(s)~ds = 0.
\]
Then
\[
	\int_{-1}^1 P_{\ell+2}(s+2)Q(s)~ds \neq 0.
\]}

{\bf Proof.} By Rolle's theorem and the well-known formula
\[
	\frac{d}{ds}P_{\ell+1}(s) = \frac{2P_\ell(s)}{\| P_\ell \|_2^2} +  \frac{2P_{\ell-2}(s)}{\| P_{\ell-2} \|_2^2} + \dots~,
\]
it is straightforward to check that for $s\in(-1,1)$ and $\ell \geq 2$ the following holds:
\begin{itemize}
\item $\displaystyle P_\ell(s+2)> 1, \quad P_\ell'(s+2)> 0, \quad P_\ell''(s+2) > 0,$
\item $\displaystyle P_{\ell + 2}(s+2)> P_\ell(s), \quad P_{\ell+2}'(s+2)> P_\ell'(s), \quad P_{\ell+2}''(s+2) > P_\ell(s+2).$
\end{itemize}
Then, if for $0\neq Q$ and for some $\ell > N$ we have
\[
	\int_{-1}^1 P_\ell(s+2)Q(s)~ds = 0,
\]
it follows that
\[
	\int_{-1}^1 P_{\ell+2}(s+2)Q(s)~ds \neq 0.
\]
\hfill $\Box$

{\bf Theorem.} {\it The matrix $B$ is invertible. Therefore, the reconstruction problem has a unique solution.}

{\bf Proof.}

We proceed by showing that the \emph{rows} of $B$ are linearly independent. Let
\[
	a_k = \left[ a_{k,N+1}, \dots, a_{k,3N+2} \right], \qquad
	c_k = \left[ c_{k,N+1}, \dots, c_{k,3N+2} \right]
\]
and assume that $\lambda_k,~\mu_k \in \mathbb{R},~k =0, \dots, N$, such that
\[
	 \sum_{k=0}^N ( \lambda_k a_k + \mu_k c_k) = 0.
\]
That means that for each $\ell = N+1, \dots, 3N+3$ we have
\[
\sum_{k=0}^N ( \lambda_k a_{k,\ell} + \mu_k c_{k,\ell})= \sum_{k=0}^N \left( \lambda_k  + \frac{1+(-1)^{k+\ell}}{2} \mu_k \right)a_{k,\ell} = 0,
\]
and so
\begin{equation}\label{invert condition}
\int_{-1}^1 P_\ell(s+2) \left\{ \sum_{k=0}^N  \left( \lambda_k  + \frac{1+(-1)^{k+\ell}}{2} \mu_k \right)P_k(s) \right\}~ds = 0.
\end{equation}
Let
\begin{align*}
	Q_0(s) = \sum_{k=0}^N \lambda_k P_k(s)+\sum_{\substack{k \leq N\\ k \text{ even}}} \mu_k P_k(s), \\
	Q_1(s) = \sum_{k=0}^N \lambda_k P_k(s)+\sum_{\substack{k \leq N\\ k \text{ odd}}}  \mu_k P_k(s), 
\end{align*}
then condition \eqref{invert condition} reads for all even $\ell$:
\[
	\int_{-1}^1 P_\ell(s+2)Q_0(s)~ds = 0,
\]
and for all odd $\ell$:
\[
	\int_{-1}^1 P_\ell(s+2)Q_1(s)~ds = 0.
\]
By the above Lemma, that means that $Q_0=Q_1 = 0$. Therefore, $\lambda_k = \mu_k = 0$ for all $k$ and thus the rows of $B$ are linearly independent. 
\hfill $\Box$

\section*{Acknowledgments}

The presented research has been financed by the European Research Council (ERC) under
the European Union's Seventh Framework Programme (FP7/2007-2013) with the research
project STiMulUs, ERC Grant agreement no. 278267.

\bibliography{references}
\bibliographystyle{plain}

\end{document}